\documentclass[11pt,a4paper]{elsarticle}
\usepackage[a4paper,margin=1in]{geometry}
\usepackage[english]{babel}
\usepackage[utf8]{inputenc}

\usepackage{algorithm} 
\usepackage[noend]{algpseudocode} 

\usepackage{xcolor}
\usepackage{cancel}
\usepackage{enumitem}

\usepackage{graphicx}

\usepackage{hyperref} 
\usepackage{natbib} 
\usepackage{amsmath}
\usepackage{amssymb}
\usepackage{amsfonts}
\usepackage[sfdefault=cmbr,OMLmathsans]{isomath}
\usepackage{mathrsfs}        
\usepackage{amssymb} 

\newtheorem{thm}{Theorem}
\newtheorem{lem}[thm]{Lemma}
\newdefinition{rmk}{Remark}
\newproof{pf}{Proof}
\newproof{pot}{Proof of Theorem \ref{thm2}}

\usepackage{lineno} 

\providecommand{\D}{\mathbb}
\newcommand{\Transpose}{^\mathrm{T}}

\newcommand{\puc}{\mathcal{Y}}
\newcommand{\intd}[1]{\,{\mathrm d}#1}
\newcommand{\MB}[1]{\boldsymbol{\mathsfit{#1}}}
\newcommand{\Nsubtext}[1]{{N_\mathrm{#1}}}
\newcommand{\N}{N}
\newcommand{\NN}{\Nsubtext{N}}
\newcommand{\NI}{\Nsubtext{N}}

\newcommand{\density}{\ensuremath{\rho}}

\newcommand{\contrast}{\chi}

\newcommand{\vek}[1]{\mathchoice{\displaystyle\boldsymbol{#1}}
	{\textstyle\boldsymbol{#1}}{\scriptstyle\boldsymbol{#1}}
	{\scriptscriptstyle\boldsymbol{#1}}}

\newcommand{\mat}[1]{\mathchoice{\displaystyle\mathbf{#1}}
	{\textstyle\mathbf{#1}}{\scriptstyle\mathbf{#1}}
	{\scriptscriptstyle\mathbf{#1}}}

\newcommand{\materialsymbol}{C}
\newcommand{\material}{\mat{\materialsymbol}}

\newcommand{\mandeld}{\ensuremath{d_{\ast}}}
\newcommand{\stress}{\ensuremath{\vek{\sigma}}}

\newcommand{\grad}{\ensuremath{\nabla}}
\newcommand{\symgrad}{\ensuremath{\grad_{\text{s}}}}
\newcommand{\symgradM}{\ensuremath{\vek{\partial}}}
\newcommand{\strain}{\ensuremath{\vek{\varepsilon}}}
\newcommand{\macrostrain}{\ensuremath{\bar{\vek{\varepsilon}}}}
\newcommand{\dispcomponent}{\ensuremath{u}}
\newcommand{\disp}{\vek{\dispcomponent}}
\newcommand{\perdisp}{\vek{\Tilde{\dispcomponent}}}
\newcommand{\Dperdisp}{\MB{\Tilde{\dispcomponent}}}
\newcommand{\Iimpulse}{\MB{I}_{\mathrm{imp}}}
\newcommand{\Icomb}{\MB{I}_{\mathrm{comb}}}

\newcommand{\testspace}{\ensuremath{V}}
\newcommand{\qweights}{\ensuremath{W}}
\newcommand{\perdispcomponent}{\Tilde{\dispcomponent}}
\newcommand{\testf}{\Tilde{\ensuremath{v}}}
\newcommand{\perdisptest}{\vek{\testf}}
\newcommand{\intvar}{\ensuremath{ \rho}}
\newcommand{\Dmat}{\MB{B}}
\newcommand{\Nmat}{\MB{N}}
\newcommand{\NQ}{\Nsubtext{Q}}
\newcommand{\parintro}[1]{\textbf{#1}}
\newcommand{\Dmaterial}{\MB{\materialsymbol}}
\newcommand{\Drmaterial}{\MB{\materialsymbol}_{\rm ref}}
\newcommand{\rmaterial}{\mat{ {\materialsymbol}_{\rm ref}}}
\newcommand{\linsysMat}{\MB{K}}
\newcommand{\RHS}{\MB{f}}

\newcommand{\FourierTrans}{\mathcal{F}}
\newcommand{\iFourierTrans}{\mathcal{F}^{-1}}
\newcommand{\precMatrix}{\MB{M}}
\newcommand{\Kref}{\MB{K}_{\text{ref}}}

\newcommand{\Green} {\MB{G}}
\newcommand{\GreenF} {\widehat{\MB{G}}}

\newcommand{\Krefinv} {\MB{K}_{\text{ref}}^{-1}}
\newcommand{\KrefF}{\widehat{\MB{K}}_{\text{ref}}}

\newcommand{\Jacobi} {\MB{J}}
\newcommand{\Jacobihalf} {\MB{J}^{1/2}}

\newcommand{\diag}{\text{diag}}
\newcommand{\norm}[1]{\left\lVert#1\right\rVert}

\DeclareMathOperator*{\argmin}{arg\,min}

\newcommand{\changeREV}[1]{\textcolor{black}{#1}}

\title{Jacobi-accelerated FFT-based solver for smooth high-contrast data}
 \graphicspath{ {./} }
\date{\today}

\begin{document}

\author[1,4]{Martin Ladecký}
\author[2]{Ivana Pultarová}
\author[3]{François Bignonnet}
\author[1,4]{Indre Jödicke}
\author[2]{Jan Zeman}
\author[1,4]{Lars Pastewka}

\address[1]{Department of Microsystems Engineering, University of Freiburg, Georges-K\"ohler-Allee 103, 79110 Freiburg, Germany}
\address[4]{Cluster of Excellence livMatS, Freiburg Center for Interactive Materials and Bioinspired Technologies, University of Freiburg, Georges-K\"ohler-Allee 105, 79110 Freiburg, Germany}

\address[2]{Faculty of Civil Engineering, Czech Technical University in Prague, Thákurova 7, 166 29 Prague 6, Czech Republic} 

\address[3]{Nantes Université, École Centrale Nantes, CNRS, GeM, UMR 6183, F-44600 Saint-Nazaire, France}

\begin{abstract}
The computational efficiency and rapid convergence of fast Fourier transform (FFT)-based solvers render them a powerful numerical tool for periodic cell problems in multiscale modeling. On regular grids, they tend to outperform traditional numerical methods. However, we show that their convergence slows down significantly when applied to microstructures with smooth, highly-contrasted coefficients.
To address this loss of performance, we introduce a Green-Jacobi preconditioner, an enhanced successor to the standard discrete Green preconditioner that preserves the quasilinear complexity, $\mathcal{O}(N \log N)$, of conventional FFT-based solvers. 
Through numerical experiments, we demonstrate the effectiveness of the Jacobi-accelerated FFT (J-FFT) solver within a linear elastic framework.  
 For problems characterized by smooth data and high material contrast, J-FFT significantly reduces the iteration count of the conjugate gradient method compared to the standard Green preconditioner.
These findings are particularly relevant for density-based topology optimization, solvers that use adaption of the grid, or nonlinear material laws, which all introduce smooth variations in the material properties that challenge conventional FFT-based solvers.
\end{abstract}
 
\maketitle

\clearpage
\section{Introduction and Motivation}\label{sec:introduction}

Fast Fourier transform (FFT)-based solvers have become a standard numerical tool for multiscale modeling of materials. 
Initially developed for homogenization of periodic microstructures~\cite{moulinec_fast_1994,moulinec_numerical_1998}, FFT-based solvers are now used for  various simulations of  heterogeneous structures on regular grids; for an overview, see  Refs.~\cite{Lucarini_2021,Schneider2021,Gierden2022}. The term ``FFT-based solver'' is broad, encompassing various formulations, discretization approaches, iterative solution methods, and a discrete Green’s operator, efficiently implemented using the FFT algorithm to accelerate the computations.

From the initial use of the Fourier basis, the FFT-based solvers have expanded to various discretization schemes, such as the finite differences~\cite{Muller1998, Willot2014, LeuteR2021} or finite elements~\cite{schneider_fft-based_2017, Leuschner2018,Ladecky2022}. These improvements mitigated discretization errors and reduced spurious oscillations that degrade the quality of local solution fields; for an overview, see Table~1~in~Ref.~\cite{Schneider2021}.
The requirement for a regular discretization grid, intrinsic in the FFT algorithm, can be relaxed through local grid adaptation techniques~\cite{ZECEVIC2022104208,ZECEVIC2023105187,BELLIS2024116658} or
by handling composite voxels either using an effective material property~\cite{KABEL2015168} or an X-FEM enrichment~\cite{Gehrig2025}.

To solve problems discretized on a regular grids, early FFT-based solvers ~\cite{moulinec_fast_1994,moulinec_numerical_1998,Eyre1999FNS} relied on fixed-point iterative schemes. Over time, a range of linear and nonlinear iterative solvers have been introduced to improve the convergence of these matrix-free methods. For an overview, see Table 2 in~Ref.~\cite{Gierden2022} or Table 4 in~Ref.~\cite{Schneider2021}. 
 
Although the state-of-the-art FFT-based solvers and their ancestors differ  in many aspects, they all use the discrete Green's operator in their core algorithms. Whether employed as a projection operator in strain-based schemes~\cite{moulinec_fast_1994,ZeVoNoMa2010AFFTH} or as a preconditioner in displacement-based schemes~\cite{schneider_fft-based_2017,Leuschner2018,Ladecky2022}, this operator plays a crucial role in improving the conditioning of the resulting system of equations. The discrete Green's operator ensures that this conditioning remains independent of the mesh size, making FFT-based solvers particularly well-suited for problems with fine discretization and a large number of degrees of freedom (DOFs).
%
Moreover, the sparse, block-diagonal structure of the discrete Green's operator in the Fourier space enables its application via the FFT, with a quasilinear complexity of $\mathcal{O}(\NI \log \NI)$, where $\NI$ denotes the total number of nodal points.

During their  30-year history, FFT-based solvers have been applied to advanced phase-field models for crystal plasticity with fatigue cracking~\cite{LUCARINI2023107670,CHEN2019167,MA2020112781}, or topology optimization~\cite{Masayoshi2025}.
%
However, for such problems typical of smooth data with high-phase contrast, the Green's operator preconditioned FFT-based solvers may exhibit slow and suboptimal convergence, as can be observed from the results presented in 
Refs.~\cite{LUCARINI2023107670,CHEN2019167,MA2020112781,Masayoshi2025}.
%
%
This limitation motivates our current research, and is demonstrated in the following simple example.

\parintro{Motivating example.}
\begin{figure}[htb!]
	\label{fig:intro_geometry} 
        		\includegraphics[width=1.0\textwidth]{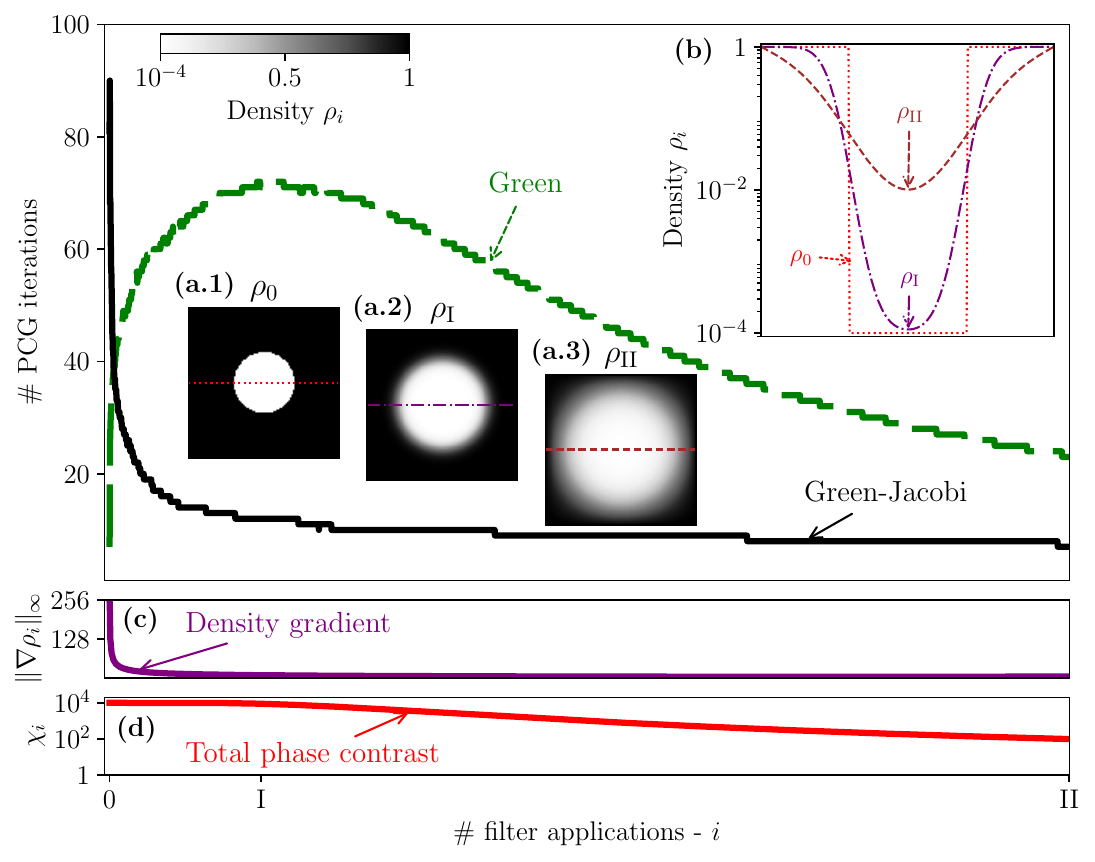} 
	\caption{Number of iterations of the preconditioned conjugate gradient (PCG) method  required to solve mechanical equilibrium 
    on a regular grid of $256^2$ nodal points as a function of the number of Gaussian filter applications,~$i$.
		The green dashed line indicates the results for Green PCG, while black solid line indicates results for Green-Jacobi PCG. 
		Figure~\textbf{(a.1)} shows the initial two-phase material density $\density_{\mathrm{0}}$, while  $\density_{\mathrm{I}}$ in~\textbf{(a.2)} is the density for which the number of iterations of Green PCG attains its maximum. The last density $\density_{\mathrm{II}}$~in~\textbf{(a.3)} has the smallest total phase contrast $\contrast_{\mathrm{II}} = 10^{2} $.
	 Panel \textbf{(b)} shows cross sections of material densities $\density_{0}$, $\density_{\mathrm{I}}$, and $\density_{\mathrm{II}}$, at the middle row of nodal points shown by the dotted, dashed-dotted, and dashed lines in panels \textbf{(a.1)} to \textbf{(a.3)}, respectively. Figure~\textbf{(c)} shows maximum gradient of density field  $\grad \density_{i}$, and Figure~\textbf{(d)} shows total phase contrast  $\max(\density_{i})/\min(\density_{i})$.  }
\end{figure}
Let us consider a compliant circular inclusion in a stiff matrix, discretized on a grid of $256^2$ nodal points. The material is linear elastic and described with a stiffness tensor $
\material(\vek{x}) 
=\density_{0}(\vek{x})\material^0,$ 
which depends on a density function $\density_{0}(\vek{x})$. The density function is $\density_{0}^\textrm{soft}=10^{-4}$ in the soft material phase and $\density_{0}^\textrm{hard}=1$ in the stiff material phase, see Figure~\ref{fig:intro_geometry}~\textbf{(a.1)}. Therefore, the initial material contrast is $\contrast_{0} = \density_{0}^\textrm{hard} / \density_{0}^\textrm{soft} = 10^{4} $. 
%
%

To construct smoother density fields $\density_{i}$, we repeatedly apply a Gaussian filter to the initial function $\density_{0}$. The filtering process involves discrete convolution of the density field with the kernel $ G=1/16\,\begin{bmatrix} 1 & 2 & 1 \end{bmatrix}^{\Transpose}\begin{bmatrix} 1 & 2 & 1 \end{bmatrix}$,
such that $\density_{i+1} = \density_{i} \star G$. The index $i =\mathrm{0}, \ldots, \mathrm{I}, \ldots ,\mathrm{II}$ indicates the number of successive filtering steps.
%

%
\changeREV{We show two-dimensional plots of two samples of filtered densities $\density_{\mathrm{I}}$ in~Figure~\ref{fig:intro_geometry}~\textbf{(a.2)}, and  $\density_{\mathrm{II}}$ in~Figure~\ref{fig:intro_geometry}~\textbf{(a.3)}. Due to the filtration, the total phase contrast decreases, as we see in~Figure~\ref{fig:intro_geometry}~\textbf{(d)} and \textbf{(b)}. In~Figure~\ref{fig:intro_geometry}~\textbf{(b)}, we show the densities $\density_{0}$, $\density_{\mathrm{I}}$, and $\density_{\mathrm{II}}$ in the middle row of nodal points. Simultaneously, the maximum of density gradient decreases, as we see in~Figure~\ref{fig:intro_geometry}~\textbf{(c)}.}  

For each density $\density_{i}$, we solve the 
micromechanical boundary value problem defined by the following system of equations:
\begin{align}
	\label{eq:intro_div_flux}
	-\grad\, \cdot \,  & \stress(\vek{x})=\vek{0}, \tag*{(mechanical equilibrium)} \notag \\
                          & \stress(\vek{x})=\material(\vek{x}) :\strain(\vek{x}), \tag*{(constitutive law)} \\
& \strain(\vek{x})=\macrostrain + \symgrad \perdisp(\vek{x}), \tag*{(kinematic compatibility)} \notag
\end{align}
where the Cauchy stress tensor $\stress(\vek{x})$ is a function of the spatially varying elastic stiffness tensors $\material(\vek{x})$ and the small strain tensor $\strain(\vek{x})$.
The small strain tensor  $\strain(\vek{x})$ is  the sum of two parts: a~constant macroscopic strain tensor $\macrostrain$ and the symmetrized gradient $\symgrad \perdisp(\vek{x})$ of the displacement fluctuation field $\perdisp(\vek{x})$, subject to periodic boundary conditions.
We use a standard continuous and piecewise linear finite element (FE) discretization on a regular grid and the preconditioned conjugate gradient (PCG) method to solve the resulting linear system. 

We compare two different preconditioners: Green, which is the standard choice in FFT-accelerated solvers, and Green-Jacobi, which is introduced and studied in the remainder of this paper. In Figure~\ref{fig:intro_geometry}, we see the number of iterations of the PCG with respect to the number of the Gauss filter applications $i$.  For the initial density $\density_{0}$ with sharp interfaces the Green PCG needs a substantially smaller number of iterations. However, for smoother densities, $\density_{i} $ with $i>0$, the number of iterations of Green PCG grows, and the Green-Jacobi method becomes a faster-converging one.
The number of iterations of the Green PCG reaches the maximum for $i=\mathrm{I}$, where the total phase contrast remains $\contrast_{\mathrm{I}} \approx 10^{4} $, i.e., close to the initial phase contrast $\contrast_0$, but the density field is smooth. Then the number of iterations of the Green PCG decreases as the total phase contrast decreases to $\contrast_{\mathrm{II}} \approx 10^{2}$.

In the article, we introduce the Green-Jacobi preconditioned FFT-based (J-FFT) solver. Through a series of numerical experiments, we illustrate problems in which the J-FFT solver outperforms the standard Green preconditioned FFT-based solver.
 
\section{Problem setup: Small-strain elasticity}\label{sec:problem_setting}
We consider a rectangular, $d$-dimensional periodic cell $\puc=\prod_{\alpha=1}^d {\left[0,l_{\alpha} \right]}$, of volume~$| \puc | =\prod_{\alpha=1}^d l_{\alpha}$, to be a~representative volume element, i.e.,~a~typical material microstructure. The symmetries of small-strain elasticity allow us to employ Mandel notation and reduce the dimension of the second-order strain tensor $\symgrad \disp=\frac{1}{2}(\grad{\disp}+\grad{\disp}^{\Transpose}): \puc \rightarrow\D{R}_{\text{sym}}^{d\times d} $ to a vector $\symgradM\disp : \puc \rightarrow \D{R}^{\mandeld}$, where $\symgradM$ is the symmetrized gradient operator, and the number of components of the symmetrized gradient in the Mandel notation is $\mandeld= (d+1)d/2 $. 
Similarly, a fourth-order tensor $ \D{C}: \puc \rightarrow\D{R}_{\text{sym}}^{d\times d\times d\times d} $ is represented by a symmetric matrix $\material : \puc \rightarrow \D{R}^{\mandeld \times \mandeld}$.
 %
%

\parintro{Strain decomposition.}
 In the small-strain micromechanical problem, the overall strain $\strain : \puc \rightarrow \D{R}^{\mandeld}$ is composed of an~average strain $\macrostrain=\frac{1}{| \puc |}\int_{\puc} \strain(\vek{x}) \intd{\vek{x}}  \in \D{R}^{\mandeld}$ and a periodically fluctuating  symmetrized gradient field $\symgradM \perdisp : \puc \rightarrow \D{R}^{\mandeld}$,
\begin{linenomath*}
\begin{align*}
    \strain(\vek{x})
    =
    \macrostrain+\symgradM \perdisp(\vek{x}) \quad \text{for all }  \vek{x} \in \puc,
\end{align*}
\end{linenomath*}
where the fluctuating displacement field $\perdisp$~belongs to the space of kinematically admissible functions $\testspace$\footnote{
$V$ is the space of vector functions $v:\mathcal{Y}\to \mathbb{R}^d$ with zero mean and such that their $\mathcal{Y}$-periodic extension 
has integrable squared derivatives on every compact 
subset of $\mathbb{R}^d.$
}.

\parintro{Weak form.} The governing equations for $ \perdisp$ are the mechanical equilibrium conditions
\begin{linenomath*}
\begin{align*}
     -\symgradM^{\Transpose}  \stress(\vek{x},\macrostrain+\symgradM \perdisp(\vek{x}))=\vek{0} \quad \text{for all }  \vek{x} \in \puc,
 \end{align*}
 \end{linenomath*}
in which $\stress:\puc \times \D{R}^{\mandeld} \rightarrow \D{R}^{\mandeld}$ is the stress field. The equilibrium equations are converted to the weak form
\begin{align}\label{eq:weak_form}
    \int_{\puc}
   \symgradM \perdisptest(\vek{x})^{\Transpose}
    \stress(\vek{x},\macrostrain+\symgradM \perdisp(\vek{x})   )
    \intd{\vek{x}}
    = 0  \quad \quad \text{for all }  \perdisptest \in \testspace,
\end{align}
where $\perdisptest$ is a test displacement field. The weak form \eqref{eq:weak_form} serves as the starting point for the discretization.
\subsection{Discretization - finite element method (FEM) on regular grid}
We discretize the weak form  \eqref{eq:weak_form}  using standard finite element method (FEM) with regular  discretization grid, as we described in~\cite{Ladecky2022}. Here  we recall the most important steps and we refer an interested reader to \cite{Ladecky2022} for more details.

\parintro{Displacement.}
 Every component $ \perdispcomponent_\alpha, \alpha=1,\dots,d,$ of the unknown displacement vector $\perdisp$ is approximated by a linear combination of  FE basis functions $\phi_{I}$.
We store the nodal values of displacement $ \perdisp(\vek{x}^{\mathrm{n}}_I)$ into a column matrix $\MB{ \Dperdisp} \in \D{R}^{d\NI} $, and write the approximation in (standard FE) matrix notation as 
\begin{align*}
        \perdispcomponent_{\alpha}(\vek{x}) \approx \perdispcomponent_{\alpha}^{N}(\vek{x}) 
        = \sum_{I=1}^{\NI}
        \phi_{I}(\vek{x}) 
        \perdispcomponent_{\alpha}^{\N}(\vek{x}^{\mathrm{n}}_{I})
        = \Nmat \MB{ \Dperdisp}_{\alpha},
    \end{align*}
where the row matrix $\Nmat: \mathcal{Y}\to \mathbb{R}^{\NI}$ stores basis functions $N_{I}=\phi_{I}(\vek{x})$, column matrix $ \MB{\Dperdisp}_{\alpha} \in \D{R}^{\NI}$ stores nodal values of the displacement in the direction ${\alpha} $, and $\NI$ denotes total number of nodal (discretization) points.
 %
 
\parintro{Strain.} Partial derivatives of this approximation are evaluated in the quadrature points~$\vek{x}^{\mathrm{q}}_Q$. The symmetrized gradient $\MB{\symgradM \Dperdisp}\in \D{R}^{\mandeld\NQ} $ at all quadrature points is given by $\symgradM \Dperdisp
=
\Dmat\MB{\Dperdisp}$,
%
where the matrix $\Dmat\in\D{R}^{\mandeld\NQ \times d\NI}$ consists of sub-matrices $\MB{B}^{\beta}\in\D{R}^{\NQ \times \NI}$ that store the partial derivatives
 \begin{linenomath*}\begin{align*}
    B^{\beta}_{Q, I}=\frac{\partial \phi_{I}}{\partial x_\beta}(\vek{x}^{\mathrm{q}}_Q
    )
     \quad \text{for } Q =1, \dots , \NQ
     \text{ and }
      I =1, \dots , \NI
     .
 \end{align*}\end{linenomath*}
Here, $\NQ$ denotes total number of quadrature points.

After the Gauss quadrature, the discretized weak form \eqref{eq:weak_form} can be rewritten in the matrix notation as 
    \begin{linenomath*}     \begin{align}\label{eq:disc_wf}
          \MB{\testf}^{\Transpose}\Dmat^{\Transpose}
           \MB{\qweights}
          \stress(\MB{E}+
         \Dmat\MB{\Dperdisp})
         =
         0
         \quad \text{for all } \MB{\testf} \in \D{R}^{d\NI}  ,
    \end{align}\end{linenomath*}
where $\MB{\testf}$ stores the nodal values of test displacements, $\MB{E} \in \D{R}^{\mandeld\NQ} $ stands for the discretized average strain, $  \stress: \D{R}^{\mandeld\NQ} \rightarrow \D{R}^{\mandeld\NQ}$ maps a vector of strains to a vector of stresses, locally at quadrature points. The diagonal matrix $\MB{\qweights}\in \D{R}^{\mandeld\NQ\times \mandeld\NQ} $ consists of $\mandeld$ identical diagonal matrices $\MB{\qweights^\text{m}}\in \D{R}^{\NQ\times \NQ}$ that store quadrature weights, $\MB{\qweights}^\text{m}_{Q,Q}=w^Q$.

Because vector $\MB{\testf}$ is arbitrary, the discretized weak form \eqref{eq:disc_wf} is equivalent to a system of  discrete nonlinear  equilibrium conditions
    \begin{align}\label{eq:non_lin_system}
          \Dmat^{\Transpose} {\MB{\qweights}}
          \stress(\MB{E}+
         \Dmat\MB{\Dperdisp})
         =
         \MB{0}.
    \end{align} 
\subsection{Linearization - Newton's method}\label{sec:linearisation}
We employ Newton's method to solve the system \eqref{eq:non_lin_system} iteratively. For this purpose, the $(i+1)$-th approximation of the nodal displacement $\MB{\Dperdisp}^{(i+1)}\in \D{R}^{d\NI}$ is given by the previous approximation $\MB{\Dperdisp}^{(i)}\in \D{R}^{d\NI}$ adjusted by a finite displacement increment~$\delta\MB{\Dperdisp}^{(i+1)}\in \D{R}^{d\NI}$,
\begin{linenomath*}\begin{align*}
    \MB{\Dperdisp}^{(i+1)}
     =
    \MB{\Dperdisp}^{(i)}+\delta\MB{\Dperdisp}^{(i+1)},
\end{align*}\end{linenomath*}
with an initial approximation $\MB{\Dperdisp}^{(0)}\in \D{R}^{d\NI}$. The displacement increment $\delta\MB{\Dperdisp}^{(i+1)}$ follows from the solution of the linear system
\begin{linenomath*}\begin{align}\label{eq:lin_system}
\underbrace{
    \Dmat^{\Transpose}
     \MB{\qweights}
    \MB{C}^{(i)} \Dmat
     }_{\linsysMat^{(i)}}
    \delta\MB{\Dperdisp}^{(i+1)}
    =
     \underbrace{
    -
    \Dmat^{\Transpose}
     \MB{\qweights}
    \stress(\MB{E}+
    \Dmat \MB{\Dperdisp}^{(i)})
     }_{\RHS^{(i)}},
\end{align}\end{linenomath*}
where the  {algorithmic} 
tangent matrix $
          \Dmaterial^{(i)}
         =
        \dfrac{\partial \stress}{\partial \strain}(\MB{E}+
         \Dmat \MB{\Dperdisp}^{(i)}) \in  \D{R}^{\mandeld\NQ\times \mandeld\NQ}, 
         $
is obtained from the constitutive tangent 
 $\material^{(i)}(\vek{x}) =\dfrac{\partial \stress}{\partial \strain}(\vek{x},\macrostrain+\symgradM \perdisp^{(i)} (\vek{x}))$, evaluated at the quadrature points  $(\vek{x}^{\mathrm{q}}_Q)$.
Traditionally, $\linsysMat^{(i)}\in  \D{R}^{d\NI\times d\NI}$  denotes the matrix of the linear system \eqref{eq:lin_system}, and $\RHS^{(i)}\in  \D{R}^{d\NI}$ stands for the right-hand side of \eqref{eq:lin_system}. 
\subsection{Linear solver - conjugate gradient (CG) method}
  For a symmetric positive-definite algorithmic tangent $\Dmaterial^{(i)}$, the system matrix $\linsysMat^{(i)} $ is symmetric and positive semi-definite, making the CG method the preferred solution method when paired with an appropriate preconditioner. 
 In the following section, we focus on an efficient preconditioning strategy for the linearized system \eqref{eq:lin_system},
 \begin{linenomath*}\begin{align*}
 		\linsysMat
 		\delta\MB{\Dperdisp}
 		=
 		\RHS,
 \end{align*}\end{linenomath*}  where
 we omit the Newton iteration index $(i)$ to improve readability.  
 \section{Preconditioning strategies}\label{sec:preconditioning}
The idea of preconditioning, see, e.g.,~\cite[Section~10.3]{golub2013matrix}
and~\cite[Chapters~9 and~10]{Saad2003}, is based on assumptions that the matrix of the preconditioned linear system
\begin{align}\label{eq:prec_lin_system_general}
\MB{M}^{-1}\linsysMat\delta\MB{\Dperdisp}=\MB{M}^{-1}\RHS,
\end{align}
has more favorable (spectral) properties 
than the original system $\linsysMat \delta\MB{\Dperdisp}=\RHS$. At the same time, the preconditioning matrix $ \precMatrix\in  \D{R}^{d\NI\times d\NI}$ should be relatively easy to invert, such that the faster convergence of the iterative method compensates for the computational overhead of the preconditioning.~\footnote{Note that system matrix $\MB{M}^{-1}\linsysMat$ is no longer symmetric. However, for symmetric $\MB{M}$ and $\linsysMat$, system~\eqref{eq:prec_lin_system_general} is equivalent with the system preconditioned in the symmetric form $\MB{M}^{-1/2}\linsysMat\MB{M}^{-1/2} \delta\MB{z}=\MB{M}^{-1/2}\RHS$, where $\delta\MB{z} =\MB{M}^{1/2}\delta\MB{\Dperdisp}$. The latter form is in fact solved when using the PCG method; see Ref.~\cite[Section~9.2.1]{Saad2003} for more details. Nonetheless, we prefer a notation with left preconditioning \eqref{eq:prec_lin_system_general} for brevity.}
 
            
\subsection{Green preconditioner}

Standard FFT-based schemes are based on a preconditioner constructed in the same way as the original matrix of the linear system~$\eqref{eq:lin_system}$,
 \begin{linenomath*} \begin{align}\label{eq:Kref}
 \Kref =
  \Dmat^{\Transpose}  \MB{\qweights}\Drmaterial  \Dmat\in  \D{R}^{d\NI\times d\NI},
 \end{align}\end{linenomath*}
where the reference algorithmic tangent matrix $\Drmaterial \in  \D{R}^{\mandeld\NQ\times \mandeld\NQ}$ corresponds to spatially uniform (constant) material data ${\rmaterial} \in  \D{R}^{\mandeld\times \mandeld}$. 

 The inverse of system matrix $ \Kref$ can be seen as a discrete Green's operator $\Green\in  \D{R}^{d\NI\times d\NI}$ of the linear system $ \Kref \, \MB{a}=\MB{b}$, i.e.,~$\Green= \Krefinv $. Notice that the spectrum of $ \Kref $ contains null eigenvalue(s) associated with rigid body translations; thus, instead of the inverse of $\Kref$, we consider its (Moore-Penrose) pseudo-inverse,\footnote{For details about the Moore-Penrose pseudo-inverse, we refer to Ref.~\cite{golub2013matrix}.} but we still denote it with $\Krefinv$ for simplicity of notation.
 

Using the discrete Green's operator $\Green$ as a preconditioner for the linear system \eqref{eq:prec_lin_system_general} leads to     
\begin{align}\label{eq:prec_lin_systemG}
\Green \linsysMat \delta\MB{\Dperdisp}=\Green\RHS,
\end{align} 
referred to as the ``Green preconditioned''.

\parintro{The fast Fourier transform.}
The system matrix $\Kref$ is block-circulant for this particular set-up involving: regular grid, spatially uniform data, and periodic boundary conditions. This implies that its discrete Fourier transform $\KrefF$ is block-diagonal and, therefore, cheap to store, cheap to multiply with, and directly, i.e cheaply, invertible in Fourier space.

Because of the above, it is common to assemble, invert, and apply the discrete Green's operator preconditioner in Fourier space using the FFT. The so-called FFT-accelerated scheme can be formally written as
\begin{linenomath*}\begin{align}\label{eq:prec_lin_systemG_FFT}
\underbrace{\iFourierTrans \GreenF \FourierTrans }_{\MB{M}^{-1}} \linsysMat \delta\MB{\Dperdisp}
=
\underbrace{\iFourierTrans \GreenF \FourierTrans }_{\MB{M}^{-1}} \RHS,
\end{align}\end{linenomath*}
where $\FourierTrans,$ and  $\iFourierTrans  $ denote the forward and inverse FFT, respectively. Multiplication with diagonal $\GreenF $ is linear in cost, $\mathcal{O}(\NI)$; therefore, the complexity of FFT $\mathcal{O}(\NI\log \NI)$ governs the overall complexity of the preconditioner. 

\changeREV{\parintro{Matrix-free assembly.} } 
\changeREV{In practice, we do not assemble matrices $\linsysMat$ and $\Kref$ explicitly. Instead, we adopt a matrix-free approach, as described in Section 5.1 of Ref.~\cite{Ladecky2022}. 
In a matrix-free implementation, we replace the system matrix with a linear operator that acts on any vector the same way as a matrix, but is computationally more efficient.
We formally replace the system matrix $ \linsysMat$ with a linear operator $\mathcal{K}: \D{R}^{d\NI} \rightarrow \D{R}^{d\NI}$, such that $\mathcal{K}  \delta\MB{\Dperdisp} =\linsysMat\delta\MB{\Dperdisp} $. Similarly, we use  a linear operator $\mathcal{K_{\text{ref}}}: \D{R}^{d\NI} \rightarrow \D{R}^{d\NI}$, such that $\mathcal{K_{\text{ref}}}  \delta\MB{\Dperdisp} =\Kref\delta\MB{\Dperdisp} $. }

\changeREV{While direct evaluation of $\GreenF$ from operator $\mathcal{K_{\text{ref}}}$  is nontrivial, it can be computed efficiently as described in Section 5.2 of Ref.~\cite{Ladecky2022}. For clarity, we present this algorithm in Algorithm~\ref{alg:green} using multi-dimensional array notation suitable for e.g. an implementation in NumPy~\cite{Numpy}.}
\changeREV{In the notation used in Algorithm~\ref{alg:green}, displacement vectors are stored as ($d+1$)-dimensional arrays of shape $[d, N_1, \dots, N_d]$, where  $   N_1 ,\dots ,N_d$ are the numbers of nodal points in the $   x_1, \dots, x_d$  directions, respectively. The algorithm consists of three steps:  i) computing the system response for all unique types of degrees of freedom, ii) computing the FFT of these responses, iii) inverting local $d \times d$ matrices in Fourier space. }
\begin{algorithm}
\caption{Assembly of the Green's Preconditioner}
\label{alg:green}
\begin{algorithmic}[1]
\State \textbf{Step 1: Compute unit impulse responses}
\State $\Kref := \mathbf{0} \in \mathbb{R}^{d \times d \times N_1 \times \dots \times N_d}$
\ForAll{components $\alpha \in \{1, \ldots, d\}$}
    \State $\Iimpulse := \mathbf{0} \in \mathbb{R}^{d \times N_1 \times \dots \times N_d}$ 
    \State $\Iimpulse[\alpha, 1, \dots, 1] = 1$ \Comment{Unit impulse}
    \State $\Kref[\alpha, :, :, \dots, :] \gets\mathcal{K_{\text{ref}}} (\Iimpulse)$ \Comment{Apply reference operator $\mathcal{K_{\text{ref}}}$}
\EndFor
\State
\State \textbf{Step 2: Fourier transform of impulse responses}
\State $\KrefF:= \mathbf{0} \in \mathbb{C}^{d \times d \times N_1 \times \dots \times N_d}$
\ForAll{$(\alpha, \beta)  \in \{0,\ldots, d\}^d$}
    \State $\KrefF[\alpha, \beta, :, \dots, :] \gets \FourierTrans\left(\KrefF [\alpha, \beta, :,\dots, :]\right)$ \Comment{Component-wise FFT}
\EndFor 
\State
\State \textbf{Step 3: Invert local matrices in Fourier space}
\State $\GreenF := \mathbf{0} \in \mathbb{C}^{d \times d \times N_1 \times \dots \times N_d}$
\ForAll{$(i_1, \dots , i_d) \in \{1, \ldots, N_1\} \times \dots  \times \{1, \ldots, N_d\}$}
    \If{$i_1 = \dots = i_d = 1$}
        \State $\GreenF[:, :, 1, 1, 1] \gets \mathbf{0}$ \Comment{Enforce zero-mean constraint}
    \Else
        \State $\GreenF[:, :, i_1, \dots ,  i_d] \gets \left(\KrefF[:, :, i_1, \dots ,  i_d]\right)^{-1}$ \Comment{Invert $d \times d$ matrices}
    \EndIf
\EndFor
\end{algorithmic}
\end{algorithm}
\changeREV{
 }

\subsection{Jacobi preconditioner}\label{sec:Jacobi}
Another basic type of preconditioner is a diagonal scaling, or the Jacobi preconditioner, which is computationally inexpensive and easy to implement, see Section~10.2~in Ref.~\cite{Saad2003}. 
 This approach is based on a preconditioner constructed from the inverse of the diagonal of the original matrix of the linear system~$\eqref{eq:lin_system}$,
 \begin{linenomath*} \begin{align}\label{eq:Jacobi}
 \Jacobi = (\diag (\linsysMat))^{-1} \in  \D{R}^{d\NI\times d\NI},
  \quad
    \Jacobi =
  \begin{bmatrix}
   \frac{1}{K_{1,1} }  & &  \\
    & \ddots & \\
    & &  \frac{1}{K_{d \NI,d \NI}}
  \end{bmatrix}.
 \end{align}\end{linenomath*}

 While Green's preconditioning is a global method, since it accounts for interactions among all DOFs across the entire domain, the Jacobi preconditioning takes into account only local interactions, making it a local method.  The Jacobi preconditioner $\Jacobi$ also incorporates information from local material data from the original problem.

 Using $\Jacobi$ from \eqref{eq:Jacobi} as a preconditioner for the linear system \eqref{eq:prec_lin_system_general} leads to the preconditioned linear system 
\begin{linenomath*}\begin{align}\label{eq:prec_lin_systemJ}
\underbrace{\Jacobi}_{\MB{M}^{-1}} \linsysMat \delta\MB{\Dperdisp}
=
\underbrace{\Jacobi}_{\MB{M}^{-1}} \RHS.
\end{align}\end{linenomath*}
The Jacobi preconditioner is diagonal in real space; therefore, its application has linear complexity $\mathcal{O}(\NI)$, and, in addition, parallelization is trivial. 

For materials with voids, where some elements of $\diag (\linsysMat)$ are zeros, we set all zero elements of $\diag (\linsysMat)$ to ones to avoid division by zeros. 
Since these DOFs correspond to void regions, they do not contribute to the solution, and the choice of replacement value does not affect convergence or accuracy. We verified numerically that varying the replacement value (e.g., using $10^{-15}$, $10^{-14}, \dots, $ $10^{15}$) does not change the number of iterations. A Jupyter notebook demonstrating this independence is included in the supplementary materials, see Ref~\cite{ladecky2025jfft}.
 
\parintro{Matrix-free assembly.}  
Direct extraction of the elements of $\diag (\linsysMat)$ from the operator $\mathcal{K}$ becomes non-trivial. A single diagonal entry can be computed via a matrix-vector product $K_{\alpha I, \alpha I} =(\mathcal{K}\MB{e}_{\alpha I} )_{\alpha I}$, where $\MB{e}_{\alpha I} \in  \D{R}^{d\NI}$ is a unit impulse vector. Unit impulse vector $\MB{e}_{\alpha I}$ has only one non-zero element equal to 1 in $\alpha$ direction in the  $I$-th nodal point. However, this requires $\NN$ matrix-vector products, so the whole process has quadratic  $ \mathcal{O}(\NN^2)$ complexity.

However, we can take advantage of the sparsity of the system matrix $\linsysMat$ (locality of the interactions/supports of basis functions) and obtain multiple diagonal terms by a single matrix-vector product. In our case, for linear finite elements, we can compute $\NI/ (d 2^d)$ terms of the diagonal at once.  
%
%
%
As a result, we assemble all elements of $\diag (\linsysMat)$ by $d 2^d$ matrix-vector products (applications of the operator $\mathcal{K}$), while maintaining the linear  $ \mathcal{O}(\NN )$ complexity.
\changeREV{The algorithm for assembling $\Jacobi$ on grids with an even number of nodal points is outlined in Algorithm~\ref{alg:jacobi}.
}
 \begin{algorithm}
\caption{Assembly of the Jacobi Preconditioner}
\label{alg:jacobi}
\begin{algorithmic}[1] 
\State \textbf{Step 1: Compute diagonal} $\linsysMat_{\diag} = \diag (\linsysMat)$
\State $\linsysMat_{\diag} := \mathbf{0} \in \mathbb{R}^{d \times N_1 \times \dots \times N_d}$ \Comment{Diagonal storage}
\State Index sets: $I_{\alpha} = \{1, \ldots, N_{\alpha}/2\}$  for $\alpha \in \{1, \ldots, d\}$
\State
\ForAll{components $\alpha \in \{1, \ldots, d\}$}
    \ForAll{offsets $(i_1, \dots , i_d) \in \{0,1\}^d$}
        \State $\Icomb := \mathbf{0} \in \mathbb{R}^{d \times N_1 \times \dots \times N_d}$ \Comment{Unit impulse comb}
        \State $\Icomb[\alpha, 2I_1-i_1, \dots , 2I_d-i_d] = 1$ \Comment{Place impulses on strided grid}
        \State $\Icomb \gets \mathcal{K}(\Icomb)$ \Comment{Apply operator $\mathcal{K}$}
        \State $\linsysMat_{\diag}[\alpha, 2I_1-i_1, \dots , 2I_d-i_d] = \Icomb[\alpha,2I_1-i_1, \dots , 2I_d-i_d]$ 
        \Comment{Extract diagonal}
    \EndFor
\EndFor
\State  \textbf{Step 2: Invert diagonal} 
\State $\Jacobi \gets \linsysMat_{\diag}^{-1}$ \quad (or $\Jacobi^{1/2} \gets \linsysMat_{\diag}^{-1/2}$ ) \Comment{Element-wise inversion}
\end{algorithmic}
\end{algorithm}

\subsection{Green-Jacobi  preconditioner} 
The last preconditioning technique examined is a synthesis of the previous two, specifically Green (global) and Jacobi (local). Integrating Green and Jacobi preconditioners can accelerate iterative solvers for extensive linear systems by combining the advantages of local preconditioners, which are cost-effective and target detailed features, with the benefits of global preconditioners, focusing on the problem's overall structure.

The findings of Gergelits et al.~\cite{Gergelits_2019}, and follow-up study~\cite{Ladecky2021} show that Green (Laplace) preconditioning yields a matrix, close to a diagonal matrix with eigenvalues equal to local material properties. Therefore, a further diagonal (Jacobi) scaling appears  to be  a meaningful strategy.  

\changeREV{The efficiency of Green-Jacobi preconditioning was theoretically estimated for certain problems involving smooth data in~Ref.~\cite[Lemma~3.2]{Serra1999}. For material data $\material(\vek{x})$ possessing two continuous derivatives, the vast majority of the eigenvalues of the Green-Jacobi preconditioned system are clustered around $1$ (or, more generally, around some positive constant).  Furthermore, the radius of this cluster increases with the amplitude of the first and second derivatives of $\material(\vek{x})$, and decreases with the amplitude of $\material(\vek{x})$ itself. This observation is consistent with the results shown in Figure~\ref{fig:intro_geometry}, where J-FFT converges faster for small norms of the density gradient. We provide a proof, inspired by Serra~Ref.~\cite{Serra1999}, in ~\ref{app:one_d_spectra} for the one-dimensional case. For more general problems with discontinuous data in two and three dimensions, we are preparing a separate publication. 
}

As we want to keep the resulting preconditioner symmetric, we split the Jacobi preconditioner $\Jacobi$ into two identical matrices
$ \Jacobi=\Jacobihalf \Jacobihalf$, where 
 \begin{linenomath*}
 \begin{align*}
\Jacobihalf 
\in  \D{R}^{d\NI\times d\NI},
  \quad
    \Jacobihalf =
  \begin{bmatrix}
    \frac{1}{\sqrt{K_{1,1}} } & &  \\
    & \ddots & \\
    & & \frac{1}{\sqrt{K_{d \NI,d \NI}} }
  \end{bmatrix}.
 \end{align*}\end{linenomath*}
Next, we wrap the Green into Jacobi preconditioners:
\begin{linenomath*}\begin{align}\label{eq:prec_lin_systemJG}
\underbrace{\Jacobihalf \Green \Jacobihalf }_{\MB{M}^{-1}} \linsysMat \delta\MB{\Dperdisp}=\underbrace{\Jacobihalf \Green \Jacobihalf }_{\MB{M}^{-1}} \RHS.
\end{align}\end{linenomath*}
%

The Jacobi preconditioner is  diagonal in the real space, while Green's preconditioner is (block) diagonal in Fourier space. Their direct combination is neither diagonal in real space nor in Fourier space. However, we can apply them sequentially in the so-called \emph{matrix-free} manner.  
The Green-Jacobi preconditioned, FFT-accelerated scheme can then be formally written as follows, 
\begin{linenomath*}\begin{align}\label{eq:prec_lin_systemJG_2}
\underbrace{\Jacobihalf \iFourierTrans \GreenF \FourierTrans \Jacobihalf }_{\MB{M}^{-1}} \linsysMat \delta\MB{\Dperdisp}
=
\underbrace{\Jacobihalf \iFourierTrans \GreenF \FourierTrans \Jacobihalf }_{\MB{M}^{-1}} \RHS.
\end{align}\end{linenomath*}
We call the resulting method a Jacobi-accelerated FFT-based (J-FFT) solver.

The overall computational overhead of the Green-Jacobi preconditioner (J-FFT) compared to the Green preconditioner (standard FFT) is 
two diagonal scalings (multiplication by $\Jacobihalf$) per iteration, and memory usage to store this diagonal of $\Jacobihalf$.  In addition, one cannot neglect the cost of assembling the Jacobi preconditioner. 
\changeREV{The assembly of $\mathrm{diag}(K)$ requires $d 2^d$ applications of the stiffness operator $\mathcal{K}$, which in 2D amounts to 8 and in 3D to 24 matrix-vector products. This is a one-time cost that is amortized over all PCG iterations. For problems where Green-Jacobi reduces the iteration count by more than $d  2^d$, the assembly cost is compensated.
}

In the following section, we perform several numerical experiments to examine and compare these three preconditioners.
%

%
%
%
%

\section{Experiments and results}\label{sec:experiments}
This paper primarily explores the application scope of the Green-Jacobi preconditioner and illustrates scenarios where Green-Jacobi PCG outperforms Green PCG and vice versa. The first two experiments are rather academic and aim to showcase the behavior of these two preconditioners on cell problems with analytically prescribed material geometries.  \changeREV{The last three examples are more applied, related to microstructure topology optimization and nonlinear elasticity, showing two of the potential applications of the Green-Jacobi preconditioner}.
To compare the performance of these preconditioners, we will use the number of iterations~$n_\mathrm{it}$ that is needed to decrease a norm of iterative error of the solution below the prescribed tolerance $\eta^{\text{CG}}$.  The  unit cell side lengths are $l_{\alpha}=1  $ in all experiments. \changeREV{All quantities in this paper are dimensionless. Simplified python implementation and Jupyter notebooks for selected examples are available at \cite{ladecky2025jfft}.}

\parintro{Constitutive law.}
 We use a linear elastic material
 \begin{align*}
 	\stress(\vek{x},  \strain(\vek{x})  ,\intvar(\vek{x}) ) 
 	=	\density(\vek{x})\material^0  \strain(\vek{x}),
 \end{align*}
where $\material^0 \in  \D{R}^{\mandeld\times \mandeld}$ is a linear elastic tensor, and $\density(\vek{x})$ is a scalar function describing  \emph{material density}.
In index notation, $
C^0_{ijkl}
=  \lambda^0 \, \delta_{ij}\, \delta_{kl}
+ \mu^0\, (\delta_{ik}\delta_{jl}+\delta_{il}\delta_{jk})$, where
$\lambda^0= 2/3 $  is the first Lamé parameter, and $  \mu^0 =  1/2$  is the shear modulus.
This choice corresponds to a material with bulk modulus $K^0= 1$.
From the discretization of the constitutive tangent 
\begin{align*}
\material( \density) 
	=	\dfrac{\partial \stress}{\partial \strain}=\density(\vek{x})\material^0,
\end{align*}
we obtain the material data matrix  $\Dmaterial(\rho)$ for the system of linear equations.
For the Green's preconditioner, we set  the reference material data ${\rmaterial =\material^0 }$.
 
\parintro{Linear system.} 
For the given linear elastic material $\material(\vek{x}, \intvar) $, the linear system \eqref{eq:lin_system} simplifies to    
\begin{linenomath*}
\begin{align}\label{eq:lin_system_exp}
\underbrace{
    \Dmat^{\Transpose} {\MB{\qweights}}
         \MB{C}(\rho) \Dmat
     }_{\linsysMat(\rho)}
   \delta\MB{\Dperdisp}
    =
     \underbrace{
    - \Dmat^{\Transpose} {\MB{\qweights}}
    \MB{C}(\rho) \MB{\macrostrain}
     }_{\RHS(\rho)}.
\end{align}\end{linenomath*} 
We consider the macroscopic gradient $\macrostrain =[1,1,1]^{\Transpose}$, unless stated otherwise.
  %

\parintro{Material data sampling.} We always consider pixel-wise constant material data. This means that $\density(\vek{x})$ is the same for all quadrature points in the pixel. To create the geometry, we sample the material distribution function $\density(\vek{x})$ in the nodal points $\vek{x}^I$, and assign this $\density(\vek{x}^I)$ to the whole pixel.  We use the shorthand $\mathcal{G}_p$ to denote the geometry with $p$ data sampling points (pixels) in each spatial direction. The total number of sampling points (pixels) is then $p^{d}$. 

\parintro{Discretization and  grid refinement.}
The computational domain $\puc$ is always discretized on a regular discretization grid with linear triangular finite elements. We use $\mathcal{T}_{n}$ to denote the discretization with  $n $ nodal  points in each spatial direction. The total number of nodal points is then $n^{d}$,  with a total number of DOFs equal to $d n^{d}$. 

\begin{figure}[htb!]     
\includegraphics[width=0.5\textwidth]{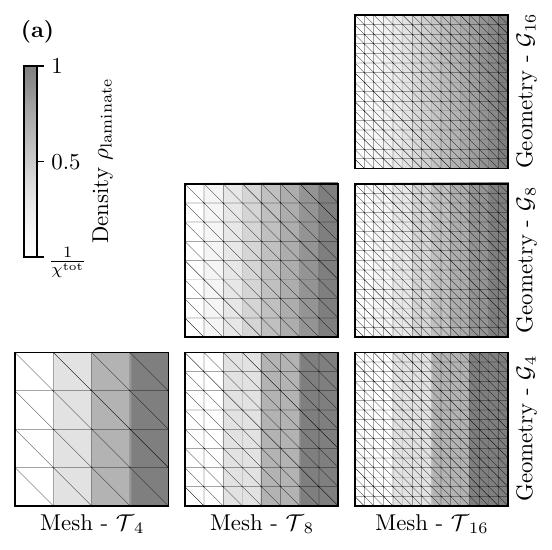} 
\includegraphics[width=0.5\textwidth]{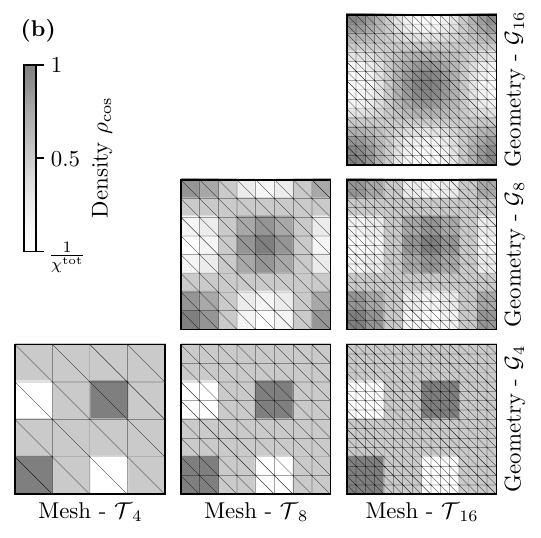} 
          \caption{Example of data samplings and discretization grids of two density functions. In \textbf{(a)}, we see \changeREV{the} linear density function $ \density_{\text{laminate}}(\vek{x})$ used in the first experiment from Section~\ref{sec:n_laminate}. In \textbf{(b)},  we see \changeREV{the} cosine density function $ \density_{\text{cos}}(\vek{x})$  used in the second experiment from Section~\ref{sec:periodic}. 
          	Both material density functions are sampled with resolutions $4^2$, $8^2$, and $16^2$ pixels denoted as  $\mathcal{G}_4$, $\mathcal{G}_8$, and $\mathcal{G}_{16}$, respectively. Finite element discretization grids  $\mathcal{T}_4$, $\mathcal{T}_8$, and $\mathcal{T}_{16}$ consist of $4^2,8^2$, and $16^2$ nodal points, respectively. The parameter $\contrast^{\rm tot}$ controls the total phase contrast}
          \label{fig:grid_refinement}
\end{figure}

The number of nodal points $n$ of $\mathcal{T}_{n}$ must be larger than or equal to the number of pixels $p$ of the geometry $\mathcal{G}_p$, that is, $n\geq p$. Otherwise, the geometry could not be captured by discretization.  An example of data sampling of two material densities with discretization grids is shown in Figure~\ref{fig:grid_refinement}.
On the left, Figure~\ref{fig:grid_refinement}  \textbf{(a)}, we see a linear function sampled with three resolutions $\mathcal{G}_4$, $\mathcal{G}_8$, and $\mathcal{G}_{16}$ with $4^2$, $8^2$, and $16^2$ pixels, respectively. These geometries are discretized on the grids $\mathcal{T}_{4}$, $\mathcal{T}_{8}$, and $\mathcal{T}_{16}$. 
On the right, Figure~\ref{fig:grid_refinement}  \textbf{(b)}, we see the same setup but for a cosine function.   

\parintro{Initial solution.}
PCG method generates a sequence of solutions $ \MB{\Dperdisp}_k$, $k=1,2,\dots$, that converges to the solution $ \MB{\Dperdisp}$ for arbitrary initial guess $ \MB{\Dperdisp}_0$. However, the total number of iterations is affected by the choice of initial guess $ \MB{\Dperdisp}_0$. 
Therefore, we set $ \MB{\Dperdisp}_0= \MB{0}$, to suppress this influence.

\parintro{Termination criteria.}
To obtain comparable results, we have to choose the appropriate termination criteria for all preconditioned schemes. That means that we have to  measure the same quantity, regardless of the preconditioning strategy. The most straightforward way is to stop the PCG iteration when the $2$-norm of the residual drops below the tolerance $ \eta^{\text{CG}}$, $\norm{\MB{r}_{k}}^2\leq \eta^{\text{CG}}$, 
 where the residual is  $\MB{r}_{k} =\RHS- \linsysMat \MB{\Dperdisp}_k$. This norm is directly accessible for all preconditioners.
 We set the tolerance $\eta^{\text{CG}}=10^{-10}$ in all  examples. 
 

 \subsection{Laminate}\label{sec:n_laminate}
 As a first example, we consider a laminate. The laminate consists of several layers of material of equal width (thickness) but different elastic properties that, in our case, depend on the material density $\density_{\text{laminate}}$. 
The material density is a linear function 
 \begin{align*}
    \density_{\text{laminate}}(\vek{x}) =  \contrast^{\rm tot}+ \frac{1-\contrast^{\rm tot}}{1-\Delta x_1}  x_1 ,
 \end{align*}
 where we call the parameter $\contrast^{\rm tot}$ the total phase contrast, and $\Delta x_1$ is the size of the geometry pixel in the $x_1$-direction. The geometry is shown in Figure~\ref{fig:grid_refinement}~\textbf{(a)}.
 
\begin{figure}[htb!]
        
                 \includegraphics[width=1.0\textwidth]{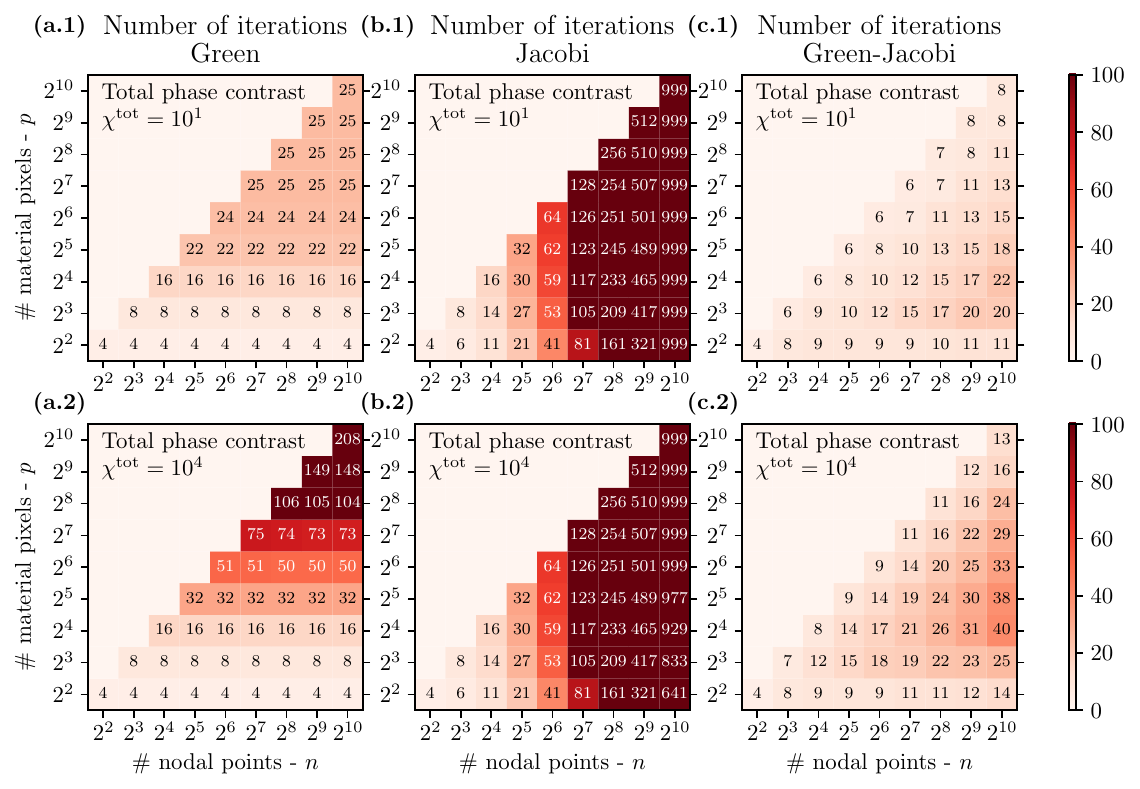}  

          \caption{Number of iterations of PCG method needed to solve mechanical equilibrium \eqref{eq:lin_system_exp}, for \changeREV{the} laminate geometry from  Section~\ref{sec:n_laminate}. Panels \textbf{(a.\_)} show results for \changeREV{the} Green preconditioner~\eqref{eq:prec_lin_systemG}, panels \textbf{(b.\_)} show results for \changeREV{the} Jacobi preconditioner~\eqref{eq:prec_lin_systemJ},  and panels \textbf{(c.\_)} show results for \changeREV{the} Green-Jacobi preconditioner~\eqref{eq:prec_lin_systemJG}. The first row \textbf{(\_.1)}  shows results for total phase contrast $\contrast^{\rm tot}=10^1$,  and  second row panels \textbf{(\_.2)} shows results for total phase contrast $\contrast^{\rm tot}=10^4$. 
          Each panel shows: i) on \changeREV{the} horizontal axis the number of nodal points in $x_1$-direction $n$ of $\mathcal{T}_{n}$, ii) on \changeREV{the} vertical axis  the number of data sampling points in $x_1$-direction $p$ of $\mathcal{G}_{p}$. \changeREV{The} upper limit for \changeREV{the} number of iterations is $999$. 
          The color coding, with color bars on the right, represents the number of iterations to highlight trends rather than exact iteration counts. }          
\label{fig:linear_iters}
\end{figure}

 \parintro{Results.}
  In Figure~\ref{fig:linear_iters}, we collect the number of iteration~$n_\mathrm{it}$ for Green  (panels \textbf{(a.\_)}), Jacobi  (panels \textbf{(b.\_)}) and Green-Jacobi  (panels \textbf{(c.\_)}) preconditioners. In the first row (panels \textbf{(\_.1)}), we see results for total phase contrasts $\contrast=10^{1}$, and in the second row  (panels \textbf{(\_.2)}) for~$\contrast=10^{4}$. In every panel, the horizontal axis indicates the number of nodal points $n$ of $\mathcal{T}_{n}$, while the vertical axis indicates the number of geometry sampling points $p$ of $\mathcal{G}_p$.  
  \begin{itemize} 
  \item  For \emph{\textbf{Green}}, first column (Figure~\ref{fig:linear_iters} - panels \textbf{(a.\_)}), the number of iterations~$n_\mathrm{it}$ remains stable as the number of nodal points $n$ of $\mathcal{T}_{n}$ increases (horizontal axis), but grows with the number of material phases (data sampling points) $p$ of $\mathcal{G}_p$ (vertical axis).
  For a small phase contrast $\contrast^{\rm tot}=10^{1}$ (Figure~\ref{fig:linear_iters} - panels \textbf{(a.1)}),  the number of iterations~$n_\mathrm{it}$ saturates at  \changeREV{$p=2^7$}   and remains stable for higher values of $p$.
  However, for the higher phase contrast $\contrast^{\rm tot}=10^{4}$ (Figure~\ref{fig:linear_iters} - panels \textbf{(a.2)}), the number of iterations continues to increase with growing $p$.
     \item For \emph{\textbf{Jacobi}}, second column (Figure~\ref{fig:linear_iters} - panels \textbf{(b.\_)}), the number of iterations~$n_\mathrm{it}$ grows with the number of nodal points $n$ of $\mathcal{T}_{n}$ (horizontal axis), but remains relatively stable with respect to the number of material phases (data sampling points) $p$ of $\mathcal{G}_p$ (vertical axis). Overall, the iteration counts for Jacobi PCG are substantially larger than those of Green PCG.
      The influence of  total phase contrast on the number of iterations~$n_\mathrm{it}$ is mild. 
               
      \item      For \emph{\textbf{Green-Jacobi}}, third column (Figure~\ref{fig:linear_iters} - panels \textbf{(c.\_s)}), we observe a relatively small increase (compared to Green) in the number of iterations~$n_\mathrm{it}$ as the number of nodal points $n$ of $\mathcal{T}_{n}$ grows (horizontal axis). In the vertical direction, the number of iterations~$n_\mathrm{it}$ decreases as the number of phases $p$ of $\mathcal{G}_p$ increases. However, the trend is not monotonous.    
      The maximum  number of iterations~$n_\mathrm{it}$ for  phase contrast $\contrast^{\rm tot}=10^{1}$ is reached for the number of material pixels  \changeREV{$p=2^4$}  discretized on a grid of $n=2^{10}$ nodal points. For phase contrast  $\contrast^{\rm tot}=10^{4}$,  we see the same  pattern.      
      In general, for higher phase contrast $\contrast^{\rm tot}=10^{4}$,  the number of iterations is higher compared to $\contrast_{\rm tot}=10^{1}$, but less than by a factor of two. 
  \end{itemize}

%
 \subsection{Cosine geometry with voids (infinite contrast).}\label{sec:periodic}
The goal of the second experiment is twofold: i) show that the absence of a large jump of the material data does not accelerate the convergence, and ii) all studied preconditioners can handle problems with infinite contrast, i.e., zero stiffness regions or voids.    

The material distribution is the cosine function 
 \begin{align*}
     \density_{\text{cos}}(\vek{x}) =  0.5 + 0.25(\cos(  2\pi (x_1-  x_2))  + \cos(  2\pi( x_2+x_1))) + \dfrac{1}{\contrast^{\rm tot}}, 
 \end{align*}
   elevated by the inverse of the total phase contrast.  
The function $ \density_{\text{cos}} $ sampled on a coarse grid~$\mathcal{G}_{4}$ results in a geometry with a gray matrix of $ \density_{\text{cos}}=0.5 +1/\contrast^{\rm tot}$, two black rigid inclusions $ \density_{\text{cos}}=1.0 +1/\contrast^{\rm tot}$, and two white soft (void) inclusions $ \density_{\text{cos}}=0.0 +1/\contrast^{\rm tot}$. The geometry is shown in Figure~\ref{fig:grid_refinement}~\textbf{(b)}. %
%
We consider two different total phase contrasts: 
First, we set $\contrast^{\rm tot}=10^{4}$ so that the density of the softest phase $\density_{\text{cos}}=10^{-4}$.  Second, we set the total phase contrast  to infinity, $\contrast^{\rm tot}=\infty $, and  therefore introduce the phase of the void with $\density_{\text{cos}}=0$. 
\begin{figure}[htb!]
                      \includegraphics[width=1.0\textwidth]{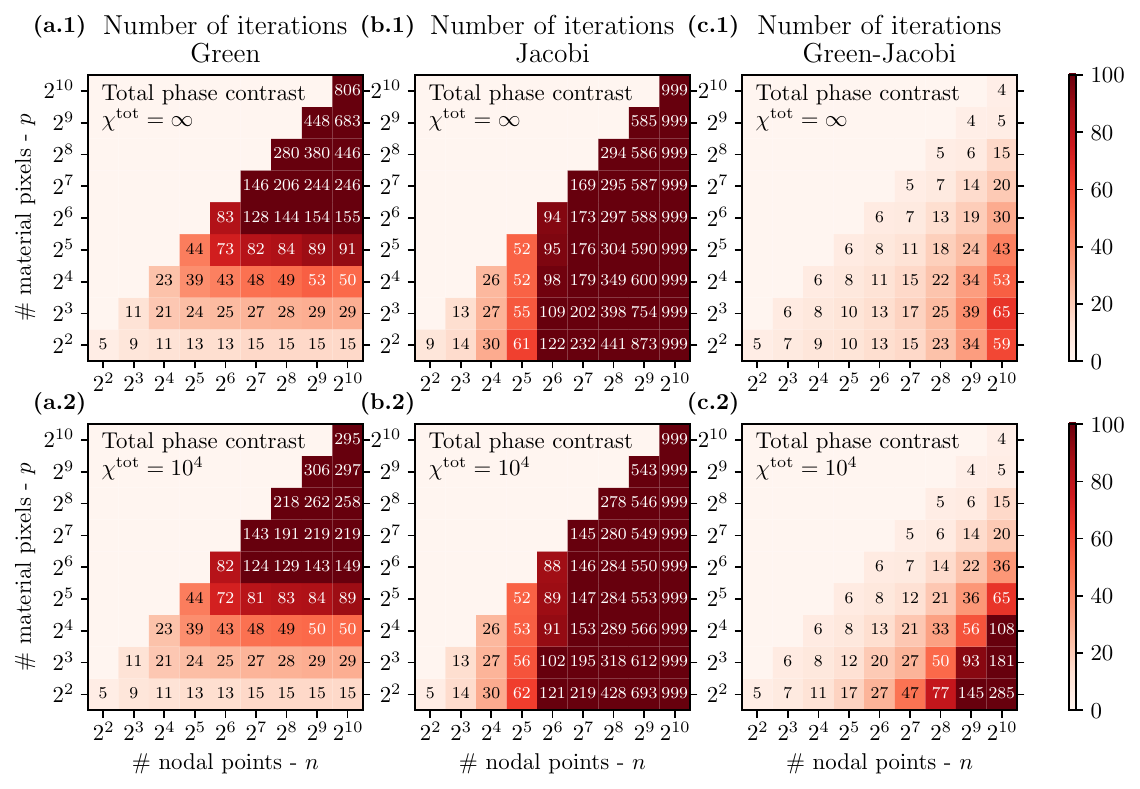}  
          \caption{
          Number of iteration of PCG method needed to solve mechanical equilibrium \eqref{eq:lin_system_exp}, for \changeREV{the} cosine geometry from  Section~\ref{sec:periodic}. Panels \textbf{(a.\_)} show results for \changeREV{the} Green preconditioner~\eqref{eq:prec_lin_systemG}, panels \textbf{(b.\_)} show results for \changeREV{the} Jacobi preconditioner~\eqref{eq:prec_lin_systemJ},  and panels \textbf{(c.\_)} show results for \changeREV{the} Green-Jacobi preconditioner~\eqref{eq:prec_lin_systemJG}. 
          	First row panels \textbf{(\_.1)} show results for total phase contrast $\contrast^{\rm tot}=\infty$,  and  second row panels \textbf{(\_.2)}  show results for total phase contrast $\contrast^{\rm tot}=10^4$. 
          Each panel has: i) on \changeREV{the} horizontal axis number of discretization points in $x_1$-direction $n$ of $\mathcal{T}_{n}$, ii) on \changeREV{the} vertical axis  the number of data sampling points in $x_1$ direction $p$ of $\mathcal{G}_{p}$.  \changeREV{The} upper limit for \changeREV{the} number of iterations is $999$.  The color coding, with color bars on the right, represents the number of iterations to highlight trends rather than exact iteration counts. }
          \label{fig:tri_iters}
\end{figure}

 \parintro{Results.}
 In Figure~\ref{fig:tri_iters}, we collect the number of iterations~$n_\mathrm{it}$ for Green (panels \textbf{(a.\_)}), Jacobi (panels \textbf{(b.\_)}) and Green-Jacobi (panels \textbf{(c.\_)}) preconditioners. In the first row (panels \textbf{(\_.1)}), we see results for total phase contrasts $\contrast^{\rm tot}=\infty$, and in the second row (panels \textbf{(\_.2)}) for $\contrast^{\rm tot}=10^{4}$.
   \begin{itemize}   
    \item 
     For \emph{\textbf{Green}}, first column (Figure~\ref{fig:tri_iters} - panels \textbf{(a.\_)}), we observe a slight increase in the number of iterations~$n_\mathrm{it}$ as the number of nodal points $n$ of $\mathcal{T}_{n}$ grows (horizontal axis). The number of iterations also increases with the number of material phases (data sampling points) $p$ of $\mathcal{G}_{p}$ (vertical axis).
     \changeREV{No notable difference is observed between materials with highly compliant inclusions and those containing voids for $p \leq 2^6$. However, for $p > 2^7$, the number of iterations required for infinite contrast begins to exceed that of finite contrast ($10^4$). At $n = 2^{10}$ and $p = 2^{10}$, the void case requires more than twice as many iterations as the finite-contrast case.} 
     \item   For \emph{\textbf{Jacobi}}, second column (Figure~\ref{fig:tri_iters} - panels \textbf{(b.\_)}), we see similar results as for \changeREV{the} laminate.  The number of iterations~$n_\mathrm{it}$ of Jacobi  increases as the number of nodal points  $n$ of $\mathcal{T}_{n}$ grows (horizontal axis), but remains relatively stable with the number of material phases (data sampling points)  $p$ of $\mathcal{G}_p$ (vertical axis). Overall, numbers of iterations for Jacobi are significantly larger than for Green.    
      \item For \emph{\textbf{Green-Jacobi}}, third  column  (Figure~\ref{fig:tri_iters} - panels \textbf{(c.\_)}), we see that the number of iterations~$n_\mathrm{it}$ increases as the number of nodal points $n$ of $\mathcal{T}_{n}$ grows (horizontal axis). However, the increase is significantly slower than for Jacobi.  In the vertical direction, when the number of phases increases, and thus the material property field becomes smoother, we see that the number of iterations~$n_\mathrm{it}$ decreases. In this experiment, these trends are monotonous, and the number of iterations~$n_\mathrm{it}$ is the highest for the coarsest geometry $\mathcal{G}_4$, and the finest discretization $\mathcal{T}_{1024}$ - (right bottom corner). Interestingly, the PCG method converges faster for infinite contrast than for $\contrast^{\rm tot}=10^4$. \changeREV{For infinite contrast, the value used to replace the inverse of the diagonal term of the stiffness matrix in the void domain has no influence on the number of iterations.} 
      

  \end{itemize}

 
%

 \subsection{Phase field topology optimization}
Microstructure topology optimization is a design approach to determine the optimal distribution of bulk material within a periodic unit cell~$\puc$. In our setting, the goal is to find a material distribution $\rho^{\text opt} $ that minimizes the difference between the target average stress response $\overline{{\stress}}_{\text{target},\gamma}$ and the macroscopic (homogenized) stress ${\stress}({x} ,\strain_{\gamma},\rho^{\text opt})$ in the microstructure, for the given set of macroscopic loads $ \macrostrain_{\gamma}$, such that  
 \begin{align}\label{eq:aim_fun}
   \rho^{\rm opt}=\argmin_{\rho} f^{\stress}(\rho)=
   \argmin_{\rho}   \sum_{\gamma}  
    \norm{\frac{1}{|\puc|}  \int_{\puc}   \;  {\stress}({x},\strain_{\gamma},\rho)     \intd{\vek{x}} - \overline{{\stress}}_{\text{target},\gamma}}^2,
\end{align}
where $\gamma$ is the load case index. For every load $ \macrostrain_{\gamma}$, we have to ensure that
the stress ${\stress}({x},\strain_{\gamma},\rho)  $ is in equilibrium. This we enforce by solving the weak mechanical equilibrium 
%
%
%
%
 in the discretized form introduced in Section~\ref{sec:problem_setting}:
\begin{linenomath*}
\begin{align}\label{eq:lin_system_topopt}
\underbrace{
    \Dmat{\MB{\qweights}}^{\Transpose}
         \MB{C}(\rho_{k}) \Dmat
     }_{\linsysMat(\rho_{k})}
    \delta\MB{\Dperdisp}
    =
     \underbrace{
    - \Dmat {\MB{\qweights}}^{\Transpose}
    \MB{C}(\rho_{k})   \MB{\macrostrain}_{\gamma}
     }_{\RHS_{\gamma}(\rho_{k})}.
\end{align}\end{linenomath*} 
Here, the index $k$ indicates an iteration number of an optimization process. We will iteratively update the density of the material $\rho_{k}$ until we reach sufficient results (a tolerance on the objective function \eqref{eq:aim_fun}).
Material data matrix  $\MB{C}(\rho_{k})$ comes again from the discretization of material data function $	\material(\density_{k}(\vek{x}))$,
similarly to the previous experiments. 

\parintro{Phase-field regularization.}
 To regularize the topology optimization problem~\eqref{eq:aim_fun}, we use the phase-field approach of ~\cite{Bourdin2003, Wallin2012}, which selects solutions with smaller interface area. We adjust the objective function \eqref{eq:aim_fun} by adding a phase-field term 
\begin{equation*}
  f^{\rm pf}(\rho, \eta) =  \eta\int_{\puc}   \;  \vert \nabla \rho \vert^2  \intd{\vek{x}} + 
  \frac{1}{\eta} \int_{\puc}   \; \rho^2 \left(1-\rho\right)^2\intd{\vek{x}}   .
\end{equation*}
The first term penalizes steep gradients in the density field $\density(\vek{x})$, promoting a smooth, constant field devoid of interphases. In contrast, the second term, represented by the double-well potential, penalizes intermediate density values ($\density(\vek{x}) $ between 0 and 1). 
As a result, the system is driven to a two-phase solution, with the density $\density(\vek{x}) \approx 0$ in the compliant material phase (representing voids), and the  $\density(\vek{x}) \approx 1$ in the stiff material phase, separated by a relatively narrow interphase.

The width of the interface is  controlled by the parameter $\eta$: A smaller $\eta$ amplifies the size of the double-well potential and therefore penalizes values of $\rho$ between $0$ and $1$, while allowing steeper gradients of $\rho$, so that the width of the diffuse interface decreases as $\eta$ decreases. 
We keep the parameter $\eta=0.01$ fixed for all resolutions.
%

\parintro{Load cases.}
We optimize for three independent load cases  $\macrostrain_\gamma = \vek{e}_\gamma$, where $\vek{e}_\gamma $ is a canonical vector satisfying  $e_{\gamma,\alpha} =\delta_{\gamma \alpha}$. 
We choose a target bulk modulus  $K_\mathrm{target}=0.025$  and target shear modulus $\mu_\mathrm{target} =  0.15$. This choice corresponds to a target material with Young modulus $E_\mathrm{target} = 0.15$ and a negative Poisson's ratio $\nu_\mathrm{target}=-0.5 $.

\parintro{Optimization algorithm.}
To minimize the objective function \eqref{eq:aim_fun}, we use gradient based optimizer, the Limited-memory Broyden–Fletcher–Goldfarb–Shanno algorithm (L-BFGS), see Ref.~\cite{LBFGS}.  
 The system is initialized with random noise,
$ \density_{0}(\vek{x})\sim \mathcal{U}(0,1)$, where the values are uniformly distributed between $0$ and $1$ without any spatial correlation. Then we iteratively update the material density function $\density_{k}(\vek{x}) $, where ${k}$ is the iteration number of L-BFGS. 
For every material distribution $\density_{k}(\vek{x}) $, we solve the micromechanical equilibrium \eqref{eq:lin_system_topopt} for each load case. We always start from the zero initial guess, i.e.,   $\MB{\Dperdisp}_0= \MB{0}$.

Topology optimization is not the main part of this article, therefore we refer the interested reader to a separated work dedicated to FFT-accelerated topology optimization~\cite{jödicke2022}. Here, our aim is to investigate the effect of a preconditioner on the convergence of PCG. Therefore, we look at the number of iterations~$n_\mathrm{it}$ of PCG needed to solve the mechanical equilibrium using the Green, Jacobi, and Green-Jacobi preconditioners. 

\begin{figure}[htb!]
\includegraphics[width=1.0\textwidth]{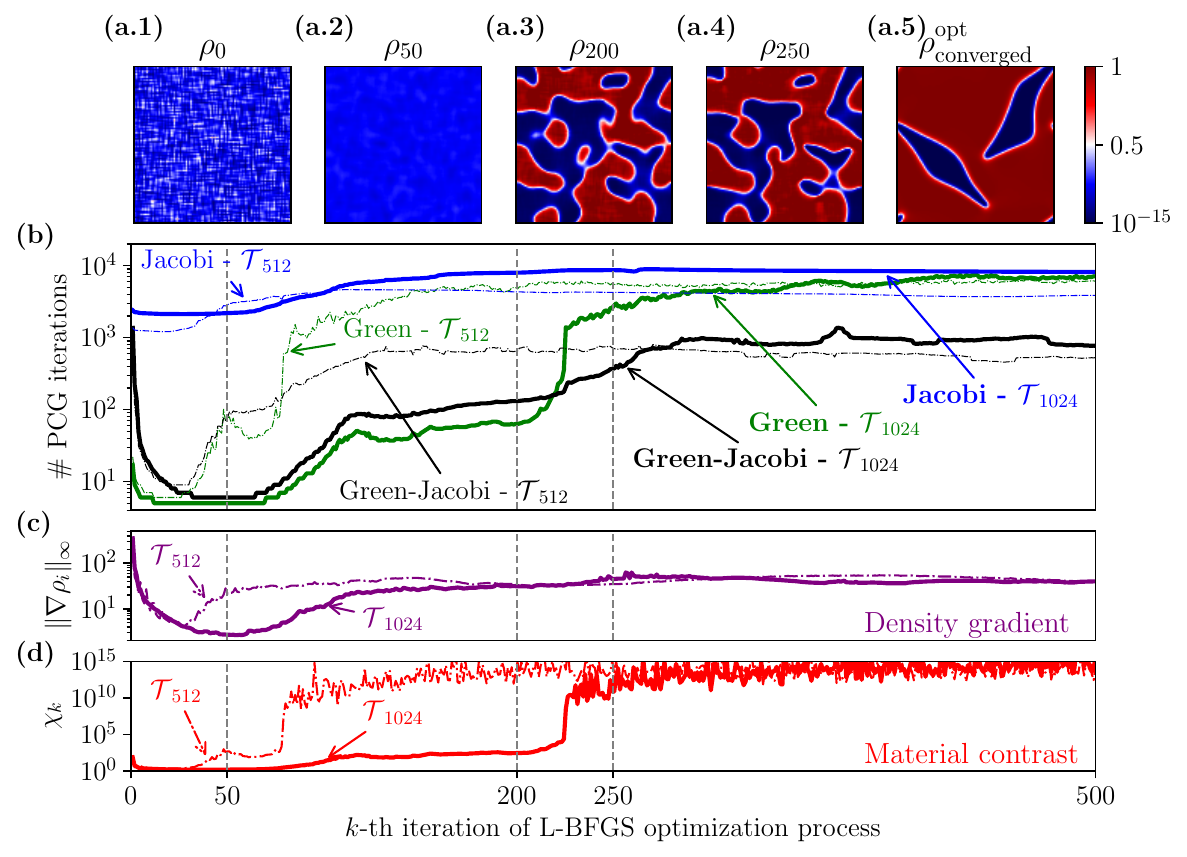}           

\caption{\changeREV{ Number of iteration of the PCG method needed to solve mechanical equilibrium \eqref{eq:lin_system_topopt}, 
with respect to the geometry obtained during the topology optimization process (using L-BFGS optimizer). In the top row we see:  initial geometry (random noise) \textbf{(a.1)}, near uniform geometry \textbf{(a.2)},  initial two-phase geometry \textbf{(a.3)}, two-phase geometry \textbf{(a.4)},  and  converged geometry \textbf{(a.5)}. Blue color indicates void (very soft material), while red indicates the bulk material. 
In the graph \textbf{(b)}, the green lines correspond to the Green preconditioner~\eqref{eq:prec_lin_systemG}, the blue lines to the Jacobi preconditioner~\eqref{eq:prec_lin_systemJ},  and  and the black lines to the Green-Jacobi preconditioner~\eqref{eq:prec_lin_systemJG}. 
Graph \textbf{(c)} displays the norm of the density gradient, while graph \textbf{(d)} shows the total phase contrast of the material.
  All three graphs shows results for 2 meshes $\mathcal{T}_{512}$, and $\mathcal{T}_{1024}$,  using dash-dotted,   and solid  lines, respectively.
  }}
\label{fig:phase_field}
\end{figure}

\parintro{Results.}
\changeREV{In~Figure~\ref{fig:phase_field}, we analyze the convergence of CG in the first $500$ steps of the L-BFGS optimization process.  After approximately $500$ iterations, the number of CG iterations stabilizes.}
\changeREV{ Figure~\ref{fig:phase_field}~\textbf{(b)} shows the number of iterations $n_\mathrm{it}$ of these three approaches, with respect to the iteration $k$ of L-BFGS optimization, for two levels of discretization $\mathcal{T}_{512},\mathcal{T}_{1024}$ with $\NI={512}^2,\NI={1024}^2$ numbers of nodes, respectively.}
\changeREV{Figure~\ref{fig:phase_field}~\textbf{(c)} shows the norm of gradient of density field $\rho_k$, and Figure~\ref{fig:phase_field}~\textbf{(d)} shows the total material phase contrast $\contrast^{\rm tot}$.}

 \changeREV{We identify four distinct stages in the evolution of the material density function  $\density_{k}(\vek{x}) $, illustrated for $\mathcal{T}_{1024}$ in  
 Figures~\ref{fig:phase_field}~\textbf{(a.\_)}. Each stage corresponds to a characteristic range of iterations. From these, we select one representative density profile per stage and annotate it for discretization $\mathcal{T}_{1024}$. Specifically we choose:} 
 \begin{itemize}   
   \item \changeREV{ $\density_{0}(\vek{x})$, corresponding to the initial random density shown in Figure~\ref{fig:phase_field}~\textbf{(a.1)}. In this configuration, the norm of gradient is large, $\|\nabla \rho_0 \|_{\infty} > 10^2 $, while the total phase contrast remains small,  $\contrast^{\rm tot} <10^2$. \emph{\textbf{Green}}  converges in fewer than $20$ iterations,  \emph{\textbf{Green-Jacobi}}  requires fewer than  $1400$, and \emph{\textbf{Jacobi}} fewer than $2400$.}
   \item \changeREV{   $\density_{50}(\vek{x})$, corresponding to the nearly spatially uniform structure shown in Figure~\ref{fig:phase_field}~\textbf{(a.2)}.   In this state, the gradient norm is small, $\|\nabla \rho_{50} \|_{\infty} < 10^1$, and the total phase contrast likewise remains low,  $\contrast^{\rm tot} <10^1$. 
     \emph{\textbf{Green}} converges in fewer than $10$ iterations,  \emph{\textbf{Green-Jacobi}} requires fewer than  $10$, and \emph{\textbf{Jacobi}}   more than $2000$}
   \item\changeREV{  $\density_{200}(\vek{x})$, corresponding to the emerging two-phase structure shown in Figure~\ref{fig:phase_field}~\textbf{(a.3)}.   In this state, the gradient norm is comparatively higher, $\|\nabla \rho_{200} \|_{\infty} < 10^2$, while the total phase contrast remains moderate,  $\contrast^{\rm tot} <10^3$. 
     \emph{\textbf{Green}} converges in fewer than $100$ iterations,  \emph{\textbf{Green-Jacobi}} requires approximately $100$, and \emph{\textbf{Jacobi}}  needs less than $8000$ iterations.}

   \item \changeREV{   $\density_{250}(\vek{x})$, corresponding to the two-phase structure  shown in Figure~\ref{fig:phase_field}~\textbf{(a.4)}.   In this state, the gradient norm is similar to the previous stage, $\|\nabla \rho_{250} \|_{\infty} < 10^2$, while the total phase contrast increases significantly, $\contrast^{\rm tot} >10^{10}$. 
     \emph{\textbf{Green}} converges in less than $3000$ iterations,  \emph{\textbf{Green-Jacobi}} requires more than $400$, and \emph{\textbf{Jacobi}}  needs more than $8000$ iterations.}

     \item \changeREV{$\rho^{\text opt}_{\text converged} $, corresponding to converged density shown in Figure~\ref{fig:phase_field}~\textbf{(a.5)}.  In this state, the gradient norm is comparable to that of the previous stage, $\|\nabla \rho_{250} \|_{\infty} < 10^2$, while the total phase contrast increases further, $\contrast^{\rm tot} >10^{12}$. For the latest stages of the optimization process, the number of iterations $n_\mathrm{it}$ of a standard FFT solver exceeds $6000$ for $\mathcal{T}_{512}$, $7500$ for $\mathcal{T}_{1024}$. The J-FFT solver needs approximately $800$ iterations for $\mathcal{T}_{512}$,   and $1400$ for $\mathcal{T}_{1024}$.  Surprisingly, the pure Jacobi preconditioner outperformed the Green preconditioner on $\mathcal{T}_{512}$, despite requiring more than $4500$ iterations. However, for $\mathcal{T}_{1024}$, Jacobi needed close to $9000$ iterations and performed worse than Green.}

 \end{itemize} 
  %
%
%

 \subsection{Smooth vs sharp interphases}\label{sec:smoothvssharp}
\begin{figure}[htb!]
\includegraphics[width=1.0\textwidth]{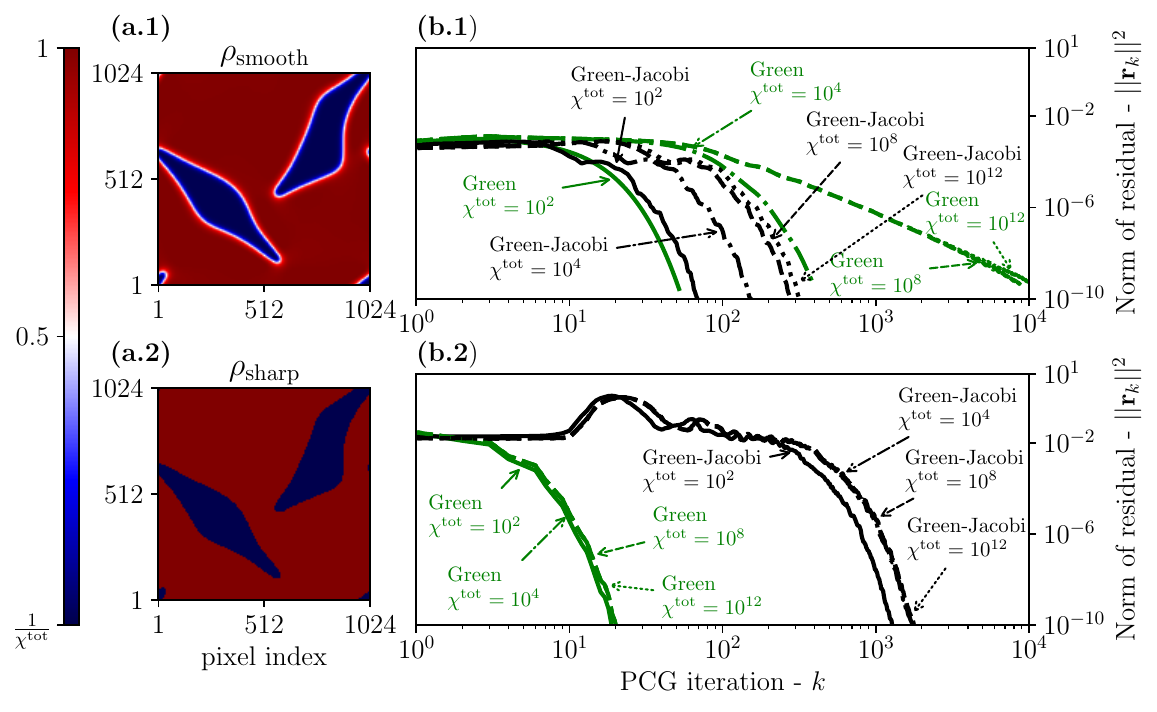}  
\caption{ Convergence of PCG for phase-field microstructure with smooth interphases  $\density_{\rm smooth}$ and sharp interphases $\density_{\rm sharp}$. In \textbf{(a.1)}, we see microstructure with smooth interphases composed of $1024^2$ pixels obtained from phase-field topology optimization.
In \textbf{(a.2)}, we see a two-phase microstructure with sharp interfaces obtained by thresholding the corresponding smooth microstructure.
	%
	 In column \textbf{(b.\_)}, we see the evolution of the   norm of residual $ \norm{\MB{r}_{k}}^2  $  with respect to the iteration number $k$ of PCG. We compare two preconditioners: i) Green (green lines) and  ii) Green-Jacobi (black lines), for total phase contrast $ \contrast^{\rm tot}=10^{2}, 10^{4}, 10^{8}$, and  $ 10^{12}$,   using solid, dash-dotted, dashed, and  dotted lines, respectively. }
\label{fig:sharp_vs_smooth_interphases}
\end{figure}

Additionally, we now compare two microstructures: one with smooth (smeared) interphases and the other with sharp interphases.
First, smooth (original) density $\density_{\rm smooth}(\vek{x})$  corresponds to the result of the phase-field topology optimization (see Figure~\ref{fig:sharp_vs_smooth_interphases}~\textbf{(a.1)}).
We rescale the  density     $\density_{\rm smooth}$ such that the total phase contrast \changeREV{$\contrast^{\rm tot}=10^{2}, 10^{4}, 10^{8}$, and  $ 10^{12}$.}
Second, a sharp density field $\density_{\rm sharp}(\vek{x})$ is obtained by thresholding of the smooth (original) density $\density_{\rm smooth}$ using the following rule, (see~Figure~\ref{fig:sharp_vs_smooth_interphases}~\textbf{(a.2)}), 
 \begin{align*}
 \density_{\rm sharp}(\vek{x})
    \begin{cases}
      1, & \text{if}\ \density_{\rm smooth}(\vek{x})  \geq  0.5, \\
     \cfrac{1}{ \contrast^{\rm tot}}, & \text{if}\ \density_{\rm smooth}(\vek{x})  <  0.5.
    \end{cases}
 \end{align*}
 Then, we solve the mechanical problems \eqref{eq:lin_system_topopt} for both $ \density_{\rm sharp}$ and $\density_{\rm smooth}$,  using Green and Green-Jacobi PCG. The evolution of the norm of the residual is shown in Figure~\ref{fig:sharp_vs_smooth_interphases}~\textbf{(b.1)} for \changeREV{a} smooth density $\density_{\rm smooth}$, and in Figure~\ref{fig:sharp_vs_smooth_interphases}~\textbf{(b.2)} for \changeREV{a} sharp density $ \density_{\rm sharp}$.  
 
 We see in Figure~\ref{fig:sharp_vs_smooth_interphases}~\textbf{(b.1)} that with increasing total phase contrast $\contrast^{\rm tot}$, the PCG convergence slows down, and the Green-Jacobi preconditioner outperforms the Green one. However, this holds only for the smooth density field $\density_{\rm smooth}$. 
In contrast, for a two-phase density field $\density_{\rm sharp}$, the convergence is almost independent of the total phase contrast $\contrast^{\rm tot}$, see Figure~\ref{fig:sharp_vs_smooth_interphases}~\textbf{(b.2)}, and the Green preconditioner now outperforms the Green-Jacobi one.
%
 \subsection{\changeREV{Nonlinear elasticity}}\label{sec:power_law}
\changeREV{
 Similar conditions, constitutive tangent with smooth variation in space, can be found in nonlinear elasticity. To demonstrate this on power-law elasticity, we adapt the three-dimensional  example from \cite[Section~4.1]{Zeman2017}.}

\begin{figure}[htb!]
\includegraphics[width=1.0\textwidth]{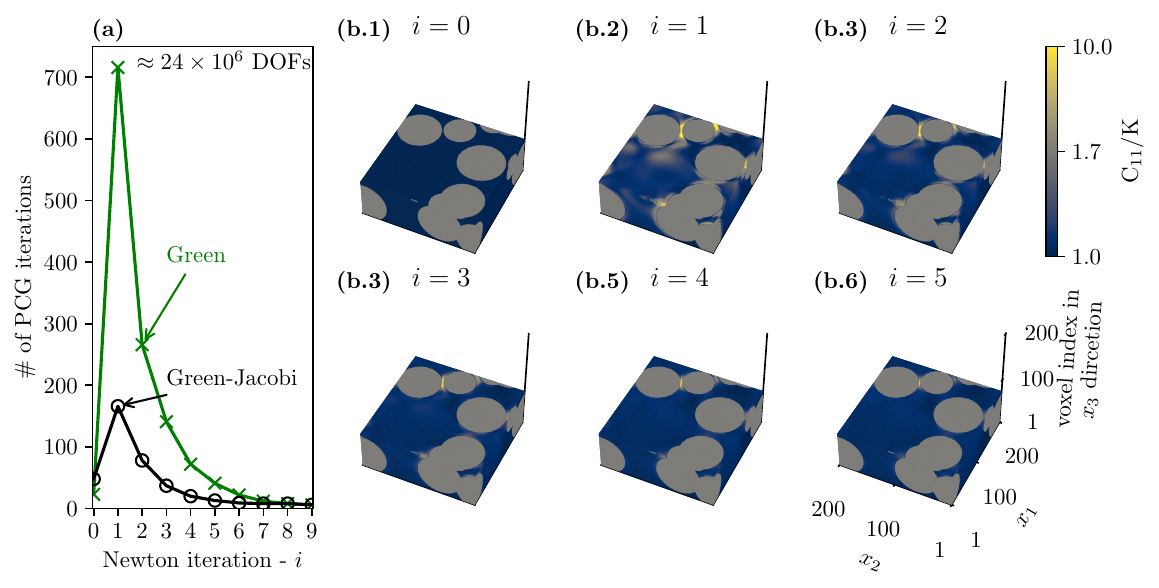}           
    
\caption{Number of   PCG iteration needed to solve linearized mechanical equilibrium \eqref{eq:lin_system}, derived for non-linear elastic constitutive problem from Section~\ref{sec:power_law}. In \textbf{(a)},  we report the number of PCG iterations as a function of the Newton iteration -  $i$.  In~\textbf{(b.\_)}, we show the $C_{11}$ component of $\material$ for the $i$-th Newton iteration, scaled by a bulk modulus $K$. For the discretization, we use trilinear finite elements with $\NI = 200^{3}$ nodal points and $8$ quadrature points per voxel/element.} 
\label{fig:nonlinear_elasticity}
\end{figure}

 \parintro{Constitutive law.}
The microstructure consists of a linear elastic spherical inclusions embedded in a non-linear matrix.
Cross-section of the geometry is shown in Figure~\ref{fig:nonlinear_elasticity}~\textbf{(b.1)}.
 The linear elastic material has the bulk modulus $K =2$, and  the shear modulus $  \mu =  1/2$. 

 For the non-linear matrix, we consider a constitutive relation
\begin{align*}
\stress= K \mathrm{tr}(\strain)  \mat{I}_{\mathrm{M}}+\sigma_\mathrm{0} \frac{2 }{3 }\left(\frac{\strain_{\mathrm{eq}}}{\strain_{\mathrm{0}}}\right)^n \frac{\strain_{\mathrm{dev}} }{\strain_{\mathrm{eq}}}  
\end{align*}
where the trace $ \mathrm{tr}(\strain)=\sum_{i=1}^{d}\strain_i$, $\mat{I}_{\mathrm{M}}$ is the second-order identity tensor in Mandel notation, 
and the equivalent strain, $ \strain_{\mathrm{eq}}$, is defined as
\begin{align*}
\strain_{\mathrm{eq}}= \sqrt{\frac{2}{3} \|\boldsymbol{\strain}_{\mathrm{dev}}\|^2 },
\end{align*}
where $\varepsilon_{\mathrm{dev}}$ is  the deviatoric strain  
\begin{align*}
\strain_{\mathrm{dev}} = \strain - \frac{1}{3}\mathrm{tr}(\strain)  \mat{I}_{\mathrm{M}}
\end{align*}
The parameters are the bulk modulus $K=2$, a reference shear stress $\sigma_\mathrm{0}=0.5$ and
strain $\varepsilon_\mathrm{0}=0.1$, and an exponent $n=5$, in analogy with \cite[Section~4.1]{Zeman2017}.
The formula for the consistent tangent operator, $ \material 
	=	\dfrac{\partial \stress}{\partial \strain}  $,  can be found in~\cite[Section~4.1]{Zeman2017}. For our purposes, it is important to note that the tangent $\material(\strain_{\mathrm{dev}}(\vek{x})) $
depends on the deviatoric strain $\strain_{\mathrm{dev}}(\vek{x})$ locally at each point $\vek{x}$.

We linearize the problem using the Newton method, as described in Section~\ref{sec:linearisation}. A macroscopic shear strain $\macrostrain$ is applied, with $\macrostrain_{6} = 0.05 \sqrt{2}$, in a single load step.
We terminated the Newton iteration once the norm of the right-hand side vector $\|\RHS^{(i)}\|$ from Eq.~\eqref{eq:lin_system} dropped below the tolerance $\eta^{\text{Newton}} = 10^{-5}$. For the discretization, we employ trilinear finite elements with $\NI = 200^{3}$ nodal points and eight quadrature points per voxel/element.

\parintro{Results.}
In Figure~\ref{fig:nonlinear_elasticity}~\textbf{(a)}, we present the number of PCG iterations for the \emph{\textbf{Green}} and \emph{\textbf{Green-Jacobi}} preconditioners at each step of the Newton method. In the first step ($i = 0$), the Green-preconditioned CG requires fewer iterations than Green-Jacobi. However, in the second step ($i = 1$), the number of iterations for Green increases significantly, exceeding that of Green-Jacobi by more than a factor of four. In the subsequent steps, the number of iterations for Green decreases and eventually reaches the same level as Green-Jacobi for $i \geq 5$.

This change in the performance of the Green-preconditioned CG arises from the following considerations. The algorithmic tangent $\material$ depends on the strain  $\strain_{\mathrm{dev}}(\vek{x})$ locally at each point~$\vek{x}$. As the strain varies smoothly across the domain, $\material$ varies accordingly. If this variation is further combined with a higher phase contrast, the efficiency of the Green-preconditioned FFT deteriorates. In Figures~\ref{fig:nonlinear_elasticity}~\textbf{(b.\_)}, we show the $C_{11}$ component of $\material$ for individual Newton iterations.

\section{Discussion}\label{sec:discusion}
In the previous section, we observed that the Green-Jacobi preconditioner outperforms the pure Green preconditioner for problems with smooth data, especially when the total phase contrast $\contrast^{\rm tot}$ is sufficiently high. Green PCG, on the other hand, outperforms the Green-Jacobi PCG for problems with a relatively small number of material phases and with sharp interphases. 

We see two main trends when using the Green preconditioner: i) the number of iterations~$n_\mathrm{it}$ is stable with respect to the mesh refinement, and ii) the number of iterations~$n_\mathrm{it}$ grows with the number of material phases~$p$ of $\mathcal{G}_{p}$. The growth of~$n_\mathrm{it}$, driven by the number of material phases~$p$ in $\mathcal{G}_{p}$, becomes increasingly significant as the total material contrast~$\contrast^{\rm tot}$ rises.
\subsection{Mesh size independence}
We experience \emph{the mesh size independence} of  \emph{\textbf{Green}} preconditioner in both academic examples: laminate from Section~\ref{sec:n_laminate} with the result in Figure~\ref{fig:linear_iters} - panels \textbf{(a.\_)} and asymptotically smooth geometry in Section~\ref{sec:periodic} with the result in Figure~\ref{fig:tri_iters} - panels \textbf{(a.\_)}. We see  that for laminate, \emph{the mesh size independence} is exactly preserved. 
Contrarily, the number of iterations~$n_\mathrm{it}$ of \emph{\textbf{Jacobi}} preconditioned system grows significantly with the system size (panels \textbf{(b.\_)}).
For the \emph{\textbf{Green-Jacobi}} preconditioner, the mesh size independence is also not preserved, as we can see in both examples (Figure~\ref{fig:linear_iters} - panels \textbf{(c.\_)} and Figure~\ref{fig:tri_iters} - panels \textbf{(c.\_)}), but the dependence is much milder than for Jacobi.
%
%
\subsection{Number of phases(interphase) dependence}
The second observed trend shows that  the number of iterations~$n_\mathrm{it}$ for \emph{\textbf{Green}} preconditioned system can be influenced by the number of material phases or interphases. Again, we can see this dependency in both academic examples: for laminate from Section~\ref{sec:n_laminate} with the result in Figure~\ref{fig:linear_iters} - panels \textbf{(a.\_)} and for geometry in Section~\ref{sec:periodic} with the result in Figure~\ref{fig:tri_iters} -  panels \textbf{(a.\_)}. We also observe that for a small phase contrast, $\contrast^{\rm{tot}}=10^1$, the number of iterations $n_\mathrm{it}$ quickly reaches a saturation point and no longer increases with the number of phases $p$ (see  Figure~\ref{fig:linear_iters} - panel \textbf{(a.1)}.

Complementary to Green PCG, the number of iterations ~$n_\mathrm{it}$ for the \emph{\textbf{Jacobi}} preconditioned system appears to be unaffected by the number of material phases~$p$ of $\mathcal{G}_{p}$.
For the \emph{\textbf{Green-Jacobi}} preconditioner, we see that the number of iterations~$n_\mathrm{it}$ decreases  when we increase the number of material phases~$p$ of $\mathcal{G}_{p}$. 
However, rather than the number of material phases $p$ in $\mathcal{G}_{p}$, we are convinced that the smoothness of the data is the key parameter  for the efficiency of the Green-Jacobi PCG.
\subsection{\changeREV{Difference between solutions.}  }\label{sec:kernels}
The solution $\MB{\Dperdisp}$ of the preconditioned linear system \eqref{eq:prec_lin_system_general} is independent of the chosen preconditioner, provided that the preconditioner shares the same kernel as the original linear system.
%
\changeREV{Preconditioning affects only the convergence rate of the iterative method, but not the exact solution of the system $\MB{\Dperdisp} $. This does not imply that the solutions $\MB{\Dperdisp}^{\mathrm{A}}_{k}$ and $\MB{\Dperdisp}^{\mathrm{B}}_{k}$, obtained with preconditioners A and B at iteration $k$, are identical; rather, they converge to the same ideal solution $\MB{\Dperdisp} $. The iterative error between these two solutions decreases as the PCG tolerance $\eta^{\text{CG}}$ approaches zero.}

\begin{figure}[htb!]     
\includegraphics[width=1.0\textwidth]{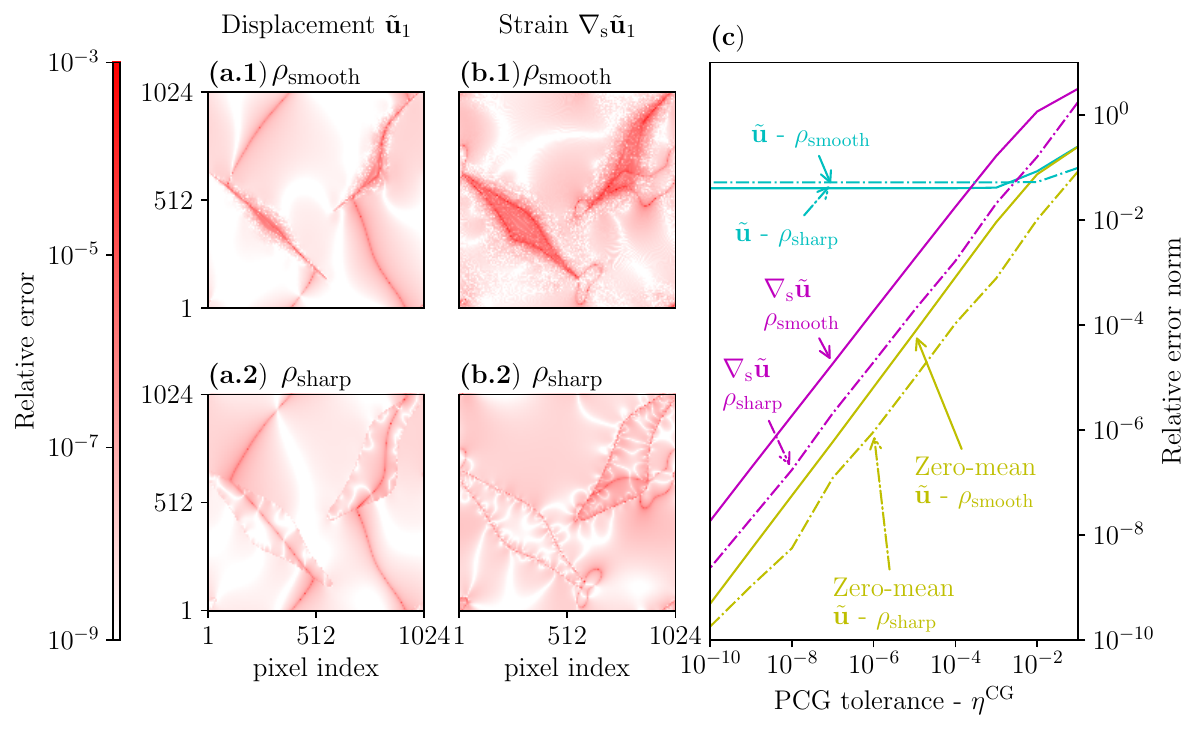}      

\caption{Relative errors between solutions obtained using the \emph{Green} PCG and the \emph{Green-Jacobi} PCG, for both sharp and smooth density fields introduced in Section~\ref{sec:smoothvssharp} with total phase contrast $\contrast^{\rm tot} = 10^{5}$.  \textbf{
(a.\_)} Relative errors in the first components of the displacement fields,  
$\left(\MB{\Dperdisp}^{\mathrm{Green-Jacobi}} - \MB{\Dperdisp}^{\mathrm{Green}}\right) / \MB{\Dperdisp}^{\mathrm{Green}}$.  
\textbf{(b.\_)} Relative errors in the first component of the symmetrized gradients,  
$\left({\symgrad \perdisp}^{\mathrm{Green-Jacobi}} - {\symgrad \perdisp}^{\mathrm{Green}}\right) / {\symgrad \perdisp}^{\mathrm{Green}}$,
for sharp and smooth density fields with total phase contrast $\contrast^{\rm tot} = 10^{5}$.
\textbf{(c)} Norm of relative error of displacements and symmetrized gradient fields as a function of the PCG tolerance $\eta^{\text{CG}}$, for both sharp and smooth density fields.  
The plots demonstrate that norm of displacement field errors stagnate at approximately $10^{-1}$ due to mean shifts, whereas subtracting the mean restores the expected convergence behavior.}
\label{fig:errors}
\end{figure}

\changeREV{To demonstrate that both the \emph{\textbf{Green}} and \emph{\textbf{Green-Jacobi}} PCGs converge to the same solution, we consider the problem introduced in Section~\ref{sec:smoothvssharp}, involving the sharp and smooth density fields $\density_{\rm sharp}$ and $\density_{\rm smooth}$. We set the total phase contrast to $\contrast^{\rm tot} = 10^{5}$.
In Figure~\ref{fig:errors}\textbf{(c)}, we present the relative error norm between the solutions  obtained using the  Green  PCG and  the Green-Jacobi  PCG, as a function of the PCG tolerance $\eta^{\text{CG}}$. Since the results for both density fields, $\density_{\rm sharp}$ and $\density_{\rm smooth}$, are similar, we discuss them jointly.}

\changeREV{As shown in Figure~\ref{fig:errors}\textbf{(c)}, the relative error norm of the symmetrized gradient $\symgrad \perdisp=  \Dmat\MB{\Dperdisp}_k$ decreases as the PCG tolerance $\eta^{\text{CG}}$ is reduced. Interestingly, the relative error norms of the displacement field $\Dperdisp=\MB{\Dperdisp}_k$ stagnate at a level of approximately $10^{-1}$. However, when the mean is subtracted from the displacement field $\MB{\Dperdisp}^{\mathrm{Green-Jacobi}}_{k}$, obtained using the  Green-Jacobi  PCG, the relative error norms of the displacement field also decrease with decreasing PCG tolerance $\eta^{\text{CG}}$, which is the expected behavior.}

\parintro{Explanation.}
In our setup, \emph{\textbf{Green}}'s preconditioner $\Green $ maps all constant fields to zero, which is fine as we require the solution $ \MB{\Dperdisp}$ to be periodic and have zero mean. This ensures, that in every iteration of Green's preconditioned CG, solution   $ \MB{\Dperdisp}^{\mathrm{Green}}_{k}$ has zero mean. 

However, this is not the case for the \emph{\textbf{Green-Jacobi}} preconditioner. Recall that the Green-Jacobi preconditioner consists of three parts: $\Jacobihalf \Green \Jacobihalf$. In the first step, the field is rescaled by $\Jacobihalf$; in the second step, it is multiplied by Green’s preconditioner $\Green$, which enforces a zero mean; and in the third step, it is rescaled again by $\Jacobihalf$, see Eq.~\eqref{eq:prec_lin_systemJG}. Since $\Jacobihalf$ is not spatially uniform over the domain, this final rescaling alters the mean of the solution field $ \MB{\Dperdisp}^{\mathrm{Green-Jacobi}}_{k}$.
 
Disregarding the iterative error, which depends on the CG tolerance $\eta^{\text{CG}}$, the two solutions are equal up to a constant vector $\MB{e}$, whose components are all identical and equal to the mean of the Green-Jacobi solution $\MB{\Dperdisp}^{\mathrm{Green-Jacobi}}_{k}$.

\parintro{Practical consequences.}  
For our micromechanical problem, it is not necessary to extract the mean from the solution, since we are primarily interested in the strain field $\MB{B}\MB{\Dperdisp}_k$. In general, the gradient field $\grad \disp$ is independent of the mean of the field $\disp$.
The same holds for the residual norm,  
$\MB{r}_{k} =\RHS- \linsysMat \MB{\Dperdisp}_k$
because the mean field lies in the kernel of the gradient operator $\MB{B}$, and therefore does not influence the magnitude of the residual norm. Nevertheless, this \textbf{ discrepancy must be kept in mind} when working directly with the displacement field.

To conclude this section, we present two-dimensional plots of the relative errors between the solutions obtained using the \emph{\textbf{Green}} PCG and the \emph{\textbf{Green-Jacobi}} PCG, for a PCG tolerance of $\eta^{\text{CG}} = 10^{-8}$.   In Figures~\ref{fig:errors}\textbf{(a.\_)} display the relative error in the first component of the displacement field,  $(\MB{\Dperdisp}^{\mathrm{Green-Jacobi}} -\MB{\Dperdisp}^{\mathrm{Green}})/\MB{\Dperdisp}^{\mathrm{Green}}$, 
while Figures~\ref{fig:errors}\textbf{(b.\_)} show the relative error in the first component of the symmetrized gradient,   $({\symgrad \perdisp}^{\mathrm{Green-Jacobi}} -{\symgrad \perdisp}^{\mathrm{Green}})/{\symgrad \perdisp}^{\mathrm{Green}}$.

\subsection{Data smoothness}  
The results indicate that the smoother material data are, the faster \emph{\textbf{Green-Jacobi}} PCG converges.
In other words, the smaller the local phase contrast is, i.e., contrast between two neighboring pixels, the faster the convergence of the Green-Jacobi PCG. 
On the other hand,  the smoother material data are, the slower \emph{\textbf{Green}} PCG converges. In our upcoming publication, we will demonstrate that it is not the smoothness, but rather the number of distinct material phases that negatively affects the convergence of the Green's PCG.

\parintro{Iterative error.}
We emphasize that our investigation focuses on the convergence of iterative solvers, namely PCG. This refers to the number of iterations required to reduce the \emph{\textbf{iterative error}}. In our setting, the iterative (or algebraic) error is orthogonal to the discretization error, which arises from discretization, as a consequence of Galerkin orthogonality. 

The results for Green PCG may seem counter-intuitive. Typically, problems with smooth data allow for a rapid reduction in discretization error compared to those with sharp discontinuities. However, in contrast to discretization error, Green PCG is unexpectedly more effective at reducing iterative error in problems with sharp interfaces between piece-wise constant data than in those with smooth data.
\section{Summary \& Conclusion}\label{sec:conclusion}
 In this paper, we discussed the efficiency of three preconditioning techniques for linear, periodic micromechanical cell problems, discretized on a regular grid using the finite element method (FEM) and solved by preconditioned conjugate gradient (PCG) method. 
  In particular, we examined two classical approaches: i) the discrete Green's operator preconditioner, which plays an essential role in FFT-based solvers, and ii) the Jacobi preconditioner, also called diagonal scaling. In addition, we introduce the J-FFT solver that employs Green-Jacobi preconditioning, which combines properties of both Green and Jacobi preconditioning techniques. 
  
  We show with numerical experiments that Jacobi is not competitive for problems with fine discretization. Green performs better for data with stronger discontinuities, while Green-Jacobi outperforms Green for smooth data with higher phase contrast.

 The computational complexity of symmetric Jacobi preconditioning is $\mathcal{O}(\NI)$, less than the complexity of Green preconditioning $\mathcal{O}(\NI \log \NI)$.
 The Green-Jacobi preconditioning requires both Green and symmetric Jacobi preconditioning, therefore it is more computationally demanding per iteration. 
However, the complexity of the J-FFT solver is still dominated by the fast Fourier transform (FFT).

Standard, discrete Green's operator preconditioned, FFT-based solvers are not well suited for problems with smoothly varying material properties. Especially, this includes problems like  
density based topology optimization~\cite{Masayoshi2025}, grid adaptation~\cite{ZECEVIC2022104208,ZECEVIC2023105187,BELLIS2024116658}, or nonlinear material laws. 
The J-FFT solver has the potential to strongly accelerate the solution of these problems.



\section*{Acknowledgments}
 ML acknowledges funding by the European Commission (Marie Skłodowska-Curie Fellowship 101106585 — microFFTTO), and Cluster of Excellence livMatS. IJ acknowledges funding by the Carl Zeiss Foundation (Research cluster ``Interactive and Programmable Materials - IPROM'') and the Deutsche Forschungsgemeinschaft (EXC 2193/1 - 390951807).  IP and JZ acknowledge funding by the European Union under the project ROBOPROX (reg. no. CZ.02.01.01/00/22\_008/0004590). 
FB thanks the staff of the Faculty of Civil Engineering at the Czech Technical University in Prague for their hospitality during his one-semester sabbatical, which was funded by Nantes Université. JZ is a member of the Ne\v{c}as Center for Mathematical Modeling.

 \section*{Declaration of generative AI and AI-assisted technologies in the writing process
}
During the preparation of this work the first author used Microsoft Copilot, and Claude in order to improve language and readability. After using this tool/service, the author reviewed and edited the content as needed and takes full responsibility for the content of the publication.

\bibliographystyle{unsrt}

\bibliography{references.bib}

@article{moulinec_fast_1994,
	author = {Moulinec, H. and Suquet, P.},
	journal = {{Comptes Rendus de l'Acad{\'e}mie des sciences. S{\'e}rie II. M{\'e}canique, physique, chimie, astronomie}},
	number = {1–2},
	pages = {1417--1423},
	title = {{A fast numerical method for computing the linear and nonlinear mechanical properties of composites}},
	volume = {318},
	year = {1994}
}

@article{moulinec_numerical_1998,
	author = {Moulinec, H. and Suquet, P.},
	doi = {10.1016/S0045-7825(97)00218-1},
	issn = {0045-7825},
	journal = {Computer Methods in Applied Mechanics and Engineering},
	number = {1–2},
	pages = {69--94},
	title = {{A numerical method for computing the overall response of nonlinear composites with complex microstructure}},
	nourl = {http://www.sciencedirect.com/science/article/pii/S0045782597002181},
	volume = {157},
	year = {1998}
}

@article{Lucarini_2021,
	author={Lucarini, S. and Upadhyay, M. V. and Segurado, J.},
	doi = {10.1088/1361-651x/ac34e1},
	year = 2021,
	publisher = {{IOP} Publishing},
	volume = {30},
	number = {2},
	pages = {023002},
	title = {{FFT} based approaches in micromechanics: {F}undamentals, methods and applications},
	journal = {Modelling and Simulation in Materials Science and Engineering},
}

@Article{Schneider2021,
	author={Schneider, M.},
	title={{A review of nonlinear {FFT}-based computational homogenization methods}},
	journal={Acta Mechanica},
	year={2021},
	day={24},
	doi={10.1007/s00707-021-02962-1},
	volume = {232},
	number = {6},
	pages = {2051--2100}
}

@Article{Gierden2022,
	author={Gierden, Ch.
	and Kochmann, J.
	and Waimann, J.
	and Svendsen, B.
	and Reese, S.},
	title={A Review of {FE-FFT}-Based Two-Scale Methods for Computational Modeling of Microstructure Evolution and Macroscopic Material Behavior},
	journal={Archives of Computational Methods in Engineering},
	year={2022},
	month={Oct},
	day={01},
	volume={29},
	number={6},
	pages={4115-4135},
	issn={1886-1784},
	doi={10.1007/s11831-022-09735-6},
}

@Article{Muller1998,
	title = "{Fourier} transforms and their application to the formation of textures and changes of morphology in solids",
	author = "Müller, {W. H.}",
	year = "1998",
	language = "English",
	pages = "61--72",
	journal= "Proc. IUTAM Symposium on Transformation Problems in Composite and Active Material"}

@article{Willot2014,
	author = {Willot, F. and Abdallah, B. and Pellegrini, Y.-P.},
	doi = {10.1002/nme.4641},
	issn = {00295981},
	journal = {International Journal for Numerical Methods in Engineering},
	month = {may},
	number = {7},
	pages = {518--533},
	publisher = {John Wiley \& Sons, Ltd},
	title = {{Fourier}-based schemes with modified {Green} operator for computing the electrical response of heterogeneous media with accurate local fields},
	nourl = {http://doi.wiley.com/10.1002/nme.4641},
	volume = {98},
	year = {2014}
}

@article{LeuteR2021,
	title = {{Elimination of ringing artifacts by finite-element projection in FFT-based homogenization}},
	journal = {Journal of Computational Physics},
	volume = {453},
	pages = {110931},
	year = {2022},
	doi = {10.1016/j.jcp.2021.110931},
	author = {R. J. Leute and {{M. Ladecký}} and A. Falsafi and I. Jödicke and I. Pultarová and J. Zeman and T. Junge and L. Pastewka}
}

@article{schneider_fft-based_2017,
	author = {Schneider, M. and Merkert, D. and Kabel, M.},
	doi = {10.1002/nme.5336},
	journal = {International Journal for Numerical Methods in Engineering},
	number = {10},
	pages = {1461--1489},
	title = {{FFT-based homogenization for microstructures discretized by linear hexahedral elements}},
	nourl = {https://onlinelibrary.wiley.com/doi/abs/10.1002/nme.5336},
	volume = {109},
	year = {2017}
}

@Article{Leuschner2018,
	author={Leuschner, M.
	and Fritzen, F.},
	title={{Fourier}-Accelerated Nodal Solvers {(FANS)} for homogenization problems},
	journal={Computational Mechanics},
	year={2018},
	month={Sep},
	day={01},
	volume={62},
	number={3},
	pages={359-392},
	issn={1432-0924},
	doi={10.1007/s00466-017-1501-5},
	url={https://doi.org/10.1007/s00466-017-1501-5}
}

@article{Ladecky2022,
	title = {An optimal preconditioned {FFT}-accelerated finite element solver for homogenization},
	journal = {Applied Mathematics and Computation},
	author = {{M. Ladecký} and  J. R. Leute and  A. Falsafi and  I. Pultarová and   L. Pastewka and  T. Junge and  J. Zeman},
	volume = {446},
	pages = {127835},
	year = {2023},
	doi = {https://doi.org/10.1016/j.amc.2023.127835}}

@article{ZECEVIC2022104208,
	title = {New large-strain {FFT}-based formulation and its application to model strain localization in nano-metallic laminates and other strongly anisotropic crystalline materials},
	journal = {Mechanics of Materials},
	volume = {166},
	pages = {104208},
	year = {2022},
	issn = {0167-6636},
	doi = {https://doi.org/10.1016/j.mechmat.2021.104208}, 
	author = {M. Zecevic and R. A. Lebensohn and L. Capolungo},
}

@article{ZECEVIC2023105187,
	title = {Non-local large-strain {FFT}-based formulation and its application to interface-dominated plasticity of nano-metallic laminates},
	journal = {Journal of the Mechanics and Physics of Solids},
	volume = {173},
	pages = {105187},
	year = {2023},
	issn = {0022-5096},
	doi = {https://doi.org/10.1016/j.jmps.2022.105187},
	author = {M. Zecevic and R. A. Lebensohn and L. Capolungo}
}

@article{BELLIS2024116658,
	title = {Numerical homogenization by an adaptive {Fourier} spectral method on non-uniform grids using optimal transport},
	journal = {Computer Methods in Applied Mechanics and Engineering},
	volume = {419},
	pages = {116658},
	year = {2024},
	issn = {0045-7825},
	doi = {https://doi.org/10.1016/j.cma.2023.116658},
	author = {C. Bellis and R. Ferrier}}

@article{KABEL2015168,
	title = {{Use of composite voxels in FFT-based homogenization}},
	journal = {Computer Methods in Applied Mechanics and Engineering},
	volume = {294},
	pages = {168-188},
	year = {2015},
	issn = {0045-7825},
	doi = {https://doi.org/10.1016/j.cma.2015.06.003},
	author = {M. Kabel and D. Merkert and M. Schneider},
}

@article{LUCARINI2023107670,
	title = {An {FFT}-based crystal plasticity phase-field model for micromechanical fatigue cracking based on the stored energy density},
	journal = {International Journal of Fatigue},
	volume = {172},
	pages = {107670},
	year = {2023},
	issn = {0142-1123},
	doi = {https://doi.org/10.1016/j.ijfatigue.2023.107670},
	author = {S. Lucarini and F. P. E. Dunne and E. Martínez-Pañeda},
}

@article{Masayoshi2025,
	author = {Matsui, M. and Hoshiba, H. and Nishiguchi, K. and Ogura, H. and Kato, J.},
	title = {Multiscale Topology Optimization Applying {FFT}-Based Homogenization},
	journal = {International Journal for Numerical Methods in Engineering},
	volume = {126},
	number = {4},
	pages = {e70009},
	doi = {https://doi.org/10.1002/nme.70009},
	year = {2025}
}

@article{CHEN2019167,
	title = {A {FFT} solver for variational phase-field modeling of brittle fracture},
	journal = {Computer Methods in Applied Mechanics and Engineering},
	volume = {349},
	pages = {167-190},
	year = {2019},
	issn = {0045-7825},
	doi = {https://doi.org/10.1016/j.cma.2019.02.017}, 
	author = {Y. Chen and D. Vasiukov and L. Gélébart and C. H. Park},
}

@article{MA2020112781,
	title = {{FFT}-based solver for higher-order and multi-phase-field fracture models applied to strongly anisotropic brittle materials},
	journal = {Computer Methods in Applied Mechanics and Engineering},
	volume = {362},
	pages = {112781},
	year = {2020},
	issn = {0045-7825},
	doi = {https://doi.org/10.1016/j.cma.2019.112781},
	author = {R. Ma and W. Sun},
}

@Article{LBFGS,
	author={Liu, D. C.
	and Nocedal, J.},
	title={On the limited memory {BFGS} method for large scale optimization},
	journal={Mathematical Programming},
	year={1989},
	day={01},
	volume={45},
	number={1},
	pages={503-528},
	issn={1436-4646},
	doi={10.1007/BF01589116},
	url={https://doi.org/10.1007/BF01589116}
}

@article{Gehrig2025,
	author = {Gehrig, F. and Schneider, M.},
	title = {An {X-FFT} Solver for Two-Dimensional Thermal Homogenization Problems},
	journal = {International Journal for Numerical Methods in Engineering},
	volume = {126},
	number = {7},
	pages = {e70022},
	doi = {https://doi.org/10.1002/nme.70022},
	year = {2025}
}

@article{Eyre1999FNS,
	author = {Eyre, D. J. and Milton, G. W.},
	journal = {The European Physical Journal Applied Physics},
	number = {1},
	pages = {41--47},
	title = {{A fast numerical scheme for computing the response of composites using grid refinement}},
	volume = {6},
	year = {1999}
}

@article{ZeVoNoMa2010AFFTH,
	author = {Zeman, J. and Vondřejc, J. and Nov{\'{a}}k, J. and Marek, I.},
	doi = {10.1016/j.jcp.2010.07.010},
	issn = {00219991},
	journal = {Journal of Computational Physics},
	number = {21},
	pages = {8065--8071},
	title = {Accelerating a {FFT}-based solver for numerical homogenization of periodic media by conjugate gradients}, 
	volume = {229},
	year = {2010}
}

@book{golub2013matrix,
	title={Matrix computations},
	author={Golub, G. H. and Van Loan, C. F.},
	isbn={9781421407944},
	lccn={2012943449},
	series={Johns Hopkins Studies in the Mathematical Sciences},
	year={2013},
	publisher={Johns Hopkins University Press}
}

@book{Saad2003,
	author = {Saad, Y.},
	title = {Iterative methods for sparse linear systems},
	publisher = {Society for Industrial and Applied Mathematics},
	year = {2003},
	doi = {10.1137/1.9780898718003},
	address = {},
	edition   = {Second},
	isbn = {978-0-898715-34-7}
}

@article{Gergelits_2019,
	author = {Gergelits, T. and Mardal, K.-A. and Nielsen, B. F. and Strakoš, Z.},
	title = {{Laplacian} Preconditioning of Elliptic {PDEs}: Localization of the Eigenvalues of the Discretized Operator},
	journal = {SIAM Journal on Numerical Analysis},
	volume = {57},
	number = {3},
	pages = {1369-1394},
	year = {2019},
	doi = {10.1137/18M1212458}
}

@article{Ladecky2021,
	author = {{{M. Ladecký}} and I. Pultarová and J.  Zeman},
	doi = {10.21136/AM.2020.0217-19},
	journal = {Applications of Mathematics},
	number = {1},
	pages = {21-42},
	publisher = {Springer Science and Business Media Deutschland GmbH},
	title = {Guaranteed Two-Sided Bounds on All Eigenvalues of Preconditioned Diffusion and Elasticity Problems Solved by the {Finite Element Method}},
	volume = {66},
	year = {2021}
}

@Article{Serra1999,
	author={Serra, S.},
	title={The rate of convergence of {Toeplitz} based {PCG} methods for second order nonlinear boundary value problems},
	journal={Numerische Mathematik},
	year={1999}, 
	day={01},
	volume={81},
	number={3},
	pages={461-495},
	issn={0945-3245},
	doi={10.1007/s002110050400}
}

@article{Bourdin2003,
	title={Design-dependent loads in topology optimization},
	author={Bourdin, B. and Chambolle, A.},
	journal={ESAIM: Control, Optimisation and Calculus of Variations},
	volume={9},
	pages={19--48},
	year={2003},
	publisher={EDP Sciences}
}

@article{Wallin2012,
	title={Optimal topologies derived from a phase-field method},
	author={Wallin, M. and Ristinmaa, M. and Askfelt, H.},
	journal={Structural and Multidisciplinary Optimization},
	volume={45},
	number={2},
	pages={171--183},
	year={2012},
	publisher={Springer}
}

@Article{jödicke2022,
	title={Efficient topology optimization using compatibility projection in micromechanical homogenization}, 
	author={I. Jödicke and R. J. Leute and T. Junge and L. Pastewka},
    	journal = {preprint}, 
	year={2022}, 
	volume = {arXiv},
		pages = {2107.04123},
	}

@article{Zeman2017,
author = {Zeman, J. and de Geus, T. W. J. and Vondřejc, J. and Peerlings, R. H. J. and Geers, M. G. D.},
title = {A finite element perspective on nonlinear {FFT}-based micromechanical simulations},
journal = {International Journal for Numerical Methods in Engineering},
volume = {111},
number = {10},
pages = {903-926}, 
doi = {https://doi.org/10.1002/nme.5481},
year = {2017}
}

@Article{Numpy,
author={Harris, C.R.
and Millman, K. J.
and van der Walt, S.J.
and Gommers, R.
and Virtanen, P.
and Cournapeau, D.
and Wieser, E.
and Taylor, J.
and Berg, S.
and Smith, N. J.
and Kern, R.
and Picus, M.
and Hoyer, S.
and van Kerkwijk, M. H.
and Brett, M.
and Haldane, A.
and del R{\'i}o, J. F.
and Wiebe, M.
and Peterson, P.
and G{\'e}rard-Marchant, P.
and Sheppard, K.
and Reddy, T.
and Weckesser, W.
and Abbasi, H.
and Gohlke, Ch.
and Oliphant, T. E.},
title={Array programming with NumPy},
journal={Nature},
year={2020},
month={Sep},
day={01},
volume={585},
number={7825},
pages={357-362}, 
issn={1476-4687},
doi={10.1038/s41586-020-2649-2}, 
}

@misc{ladecky2025jfft,
  author       = {Ladeck\'{y}, Martin and Pultarov\'{a}, Ivana and Bignonnet, Fran\c{c}ois and J\"{o}dicke, Indre and Zeman, Jan and Pastewka, Lars},
  title        = {Supplementary Code - {J-FFT} solver for smooth high-contrast data}, 
  howpublished = {\url{https://github.com/MartinLadecky/J-FFT-solver-for-smooth-high-contrast-data}}
}

\appendix

\section{\changeREV{Spectrum of the Green-Jacobi preconditioned system in 1D}}
\label{app:one_d_spectra} 

The efficiency of the Green–Jacobi preconditioner is analyzed by estimating the spectra of the resulting preconditioned matrices. In this section, we provide a detailed proof for the one-dimensional second-order differential equation. Our arguments are inspired by Serra’s work~\cite{Serra1999}, but are developed here in more detail. We assume that the material properties are described by functions with smooth second derivatives. While the approach is applicable to various types of boundary conditions, we restrict our attention to the case of periodic boundary conditions.

Let $\puc= {\left[0,1\right]}$ and let us solve 
\begin{align}
	\label{eq:Auf}
	-\dfrac{ \partial}{\partial x} \left(c(x)  \dfrac{ \partial}{\partial x} u(x)\right) =f(x)
\end{align}
where the material data function $c\in C^2(\overline{\puc})$, $c(x)\ge \delta>0$ for $ x\in\puc$,
and let $f\in L^2(\mathcal{Y})$.

We then discretize the equation using the finite element method (FEM) on a uniform grid with piecewise linear basis functions, employing $\NI$ elements of length $h = 1/\NI$. 
Let the nodal points be denoted by $x^{\mathrm{n}}_{i}$, where the index $i = 1,2,\dots,\NI$ is taken to be $\NI$-periodic. Owing to the periodic boundary conditions, we have $u(x^{\mathrm{n}}_{1}) = u(x^{\mathrm{n}}_{\NI+1})$.  
The quadrature points ${x}^{\mathrm{q}}_{i}$, for $i = 1,2,\dots,\NI$, are located at the centers of the elements  $(x^{\mathrm{n}}_{i}, x^{\mathrm{n}}_{i+1})$.

Then we obtain a system of $\NI$  linear equations
\begin{equation}\nonumber
\linsysMat\MB{u}=\RHS
\end{equation}
where $\linsysMat\in  \D{R}^{\NI\times \NI}$ is a $h$-scaled symmetric positive semi-definite matrix. 
Let the elements of $\linsysMat$ be computed using the Gauss quadrature by  
\begin{equation}\nonumber
\linsysMat_{i-1,i}
=h
\int_{x^{\mathrm{n}}_{i-1}}^{{{x}^{\mathrm{n}}_{i}}}
\dfrac{ \partial\phi_{i}(x)}{\partial x}
c(x) \,
\dfrac{ \partial\phi_{i-1}(x)}{\partial x}
 \,{\rm d}x
=
 h^2
 c(x^{\mathrm{q}}_{i-1})
\dfrac{ \partial\phi_{i}(x^{\mathrm{q}}_{i-1})}{\partial x}
 \,
\dfrac{ \partial\phi_{i-1}(x^{\mathrm{q}}_{i-1})}{\partial x} 
:=
 -1a(x^{\mathrm{q}}_{i-1})
\end{equation} 

and
\begin{align*} 
\linsysMat_{i,i}
=&
 h\int_{x^{\mathrm{n}}_{i-1}}^{{{x}^{\mathrm{n}}_{i}}}
\dfrac{ \partial\phi_{i}(x)}{\partial x}
c(x) \,
\dfrac{ \partial\phi_{i}(x)}{\partial x}
 \,{\rm d}x
 +
 h\int_{x^{\mathrm{n}}_{i}}^{{{x}^{\mathrm{n}}_{i+1}}}
\dfrac{ \partial\phi_{i}(x)}{\partial x}
c(x) \,
\dfrac{ \partial\phi_{i}(x)}{\partial x}
 \,{\rm d}x
\\
=&
 h^2
 c(x^{\mathrm{q}}_{i-1})
\dfrac{ \partial\phi_{i}(x^{\mathrm{q}}_{i-1})}{\partial x}
 \,
\dfrac{ \partial\phi_{i-1}(x^{\mathrm{q}}_{i-1})}{\partial x} 
+
 h^2
 c(x^{\mathrm{q}}_{i-1})
\dfrac{ \partial\phi_{i}(x^{\mathrm{q}}_{i-1})}{\partial x}
 \,
\dfrac{ \partial\phi_{i-1}(x^{\mathrm{q}}_{i-1})}{\partial x} 
\\
:=&
 1 a(x^{\mathrm{q}}_{i-1})
 +
 1 a(x^{\mathrm{q}}_{i})
\end{align*} 
i.e.~$\linsysMat$ is an $h$ scaled typical FEM matrix of~\eqref{eq:Auf}.
Let $\linsysMat_0$ be obtained in the same way as $\linsysMat$, but with $c(x)=1$, $x\in\mathcal{Y}$.
Note that $\Green=\linsysMat_0^{-1}$ and that 
the discretization stencils of $\linsysMat$ and $\linsysMat_0$ are  given by 
$(-a(x^{\mathrm{q}}_{i-1}), a(x^{\mathrm{q}}_{i-1})+a(x^{\mathrm{q}}_{i}),-a(x^{\mathrm{q}}_{i}))$ and
$(-1,2,-1)$, respectively.
Let $\Jacobi$ be a diagonal matrix with elements
$\Jacobi_{ii}=\frac{(\linsysMat_0)_{ii}}{\linsysMat_{ii}}$.

We aim to estimate the spectrum of the system matrix $\Jacobihalf\Green\Jacobihalf\linsysMat$
of~\eqref{eq:prec_lin_systemJG} 
which is the same as the spectrum of 
$\linsysMat_0^{-1}\Jacobihalf\linsysMat\Jacobihalf$. Let us denote $\linsysMat_{\Jacobi}=\Jacobihalf\linsysMat\Jacobihalf$.
Let $\lfloor m\rfloor$ denote the largest
integer smaller than or equal to $m$.
\begin{lem}\label{lemma1}
Let $c$, $\linsysMat$, $\Jacobi$,
$\linsysMat_0$ and $\linsysMat_{\Jacobi}$ be defined as in the previous paragraph and
let $\alpha\in(0,1)$.
Then
for every $n\in\mathbb{N}$, there exists a subspace $V_{n}\subset\mathbb{R}^{\NI}$ of dimension ${\rm dim}(V_{n})=\NI-{2\lfloor {\NI}^\alpha/2\rfloor}-1$ such that for all $\MB{v}\in V_{n}$
\begin{equation}\nonumber
\frac{\MB{v}^T
\linsysMat_{\Jacobi}
\MB{v}}{\MB{v}^T\linsysMat_0\MB{v}}= 1+O(h^{2\alpha})+O(h).
\end{equation}
The coefficient of the second  term does not depend on $\alpha$.
\end{lem}
\begin{pf}
Since $c\in C^2(\overline{\mathcal{Y}})$, then from Taylor's expansion for all $x\in\mathcal{Y}$ there exists $\xi\in(x,x+h)$ such that 
$ c(x+h)=c(x)+c'(x)h+\frac{h^2}{2}c''(\xi)$.
Therefore, $c(x^{\mathrm{q}}_{i-2})=c(x^{\mathrm{q}}_{i-1})-h c'(x^{\mathrm{q}}_{i-1})+\frac{h^2}{2}c''(\xi_{i-1}) $ with $\xi_{i-1} \in(x^{\mathrm{q}}_{i-2},x^{\mathrm{q}}_{i-1})$,    $c(x^{\mathrm{q}}_{i})=c(x^{\mathrm{q}}_{i-1})+h c'(x^{\mathrm{q}}_{i-1})+\frac{h^2}{2}c''(\zeta_{i-1}) $ with $ \zeta_{i-1}\in(x^{\mathrm{q}}_{i-1},x^{\mathrm{q}}_{i})$,
and the off-diagonal term  $(\linsysMat_{\Jacobi})_{i-1,i}=(\Jacobihalf)_{i-1,{i-1}} (\linsysMat)_{i-1,i}(\Jacobihalf)_{i,i}$.

Then,  
\begin{eqnarray}\nonumber
(\linsysMat_{\Jacobi})_{i-1,i}&=&\frac{-2c(x^{\mathrm{q}}_{i-1})}
{\sqrt{(c(x^{\mathrm{q}}_{i-2})+c(x^{\mathrm{q}}_{i-1}))(c(x^{\mathrm{q}}_{i-1})+c(x^{\mathrm{q}}_{i-1}}))} \nonumber\\
&=&
\frac{-1}{\sqrt{\left(1-\frac{hc'(x^{\mathrm{q}}_{i-1})}{2c(x^{\mathrm{q}}_{i-1})}
+\frac{h^2c''(\xi_{i-1})}{4c(x^{\mathrm{q}}_{i-1})}\right)\left(1+\frac{hc'(x^{\mathrm{q}}_{i-1})}{2c(x^{\mathrm{q}}_{i-1})}+\frac{h^2c''(\zeta_{i-1})}{4c(x^{\mathrm{q}}_{i-1})}\right)}}
\nonumber\\
&=&\frac{-1}{\sqrt{1-\frac{h^2c'(x^{\mathrm{q}}_{i-1})^2}{4c^2(x^{\mathrm{q}}_{i-1})}+
\frac{h^2c''(\xi_{i-1})}{4c(x^{\mathrm{q}}_{i-1})}+
\frac{h^2c''(\zeta_{i-1})}{4c(x^{\mathrm{q}}_{i-1})}
+O(h^3)}}\nonumber\\
&=&-1+\frac{h^2}{8}\left(\frac{c'(x^{\mathrm{q}}_{i-1})^2}{c^2(x^{\mathrm{q}}_{i-1})}
-\frac{c''(\xi_{i-1})}{c(x^{\mathrm{q}}_{i-1})}-\frac{c''(\zeta_{i-1})}{c(x^{\mathrm{q}}_{i-1})}
\right)+O(h^3),
\nonumber
\end{eqnarray}
where $\xi_{i-1},\zeta_{i-1}\in(x^{\mathrm{q}}_{i-2},x^{\mathrm{q}}_{i})$.
Thus we obtain the expansion of $\linsysMat_{\Jacobi}$ as
\begin{equation}\nonumber
\linsysMat_{\Jacobi}=\linsysMat_0+h^2\MB{T}+
O(h^3)\MB{R},
\end{equation}
where $\MB{T}$ and $\MB{R}$ are 
two-diagonal matrices with symmetric non-zero structure and with $O(1)$ elements.
Especially,  
\begin{equation}\label{Tcc}
\vert\MB{T}_{i-1,i}\vert\le  \frac{\widetilde{c}_1^2}{8\widetilde{c}_0^2}+\frac{\widetilde{c}_2}{4\widetilde{c}_0},  
\end{equation}
where $\widetilde{c}_0$ is minimum of $a$, and $\widetilde{c}_1$ and $\widetilde{c}_2$ are maxima of $\vert c'\vert$ and $\vert c''\vert$ on $\mathcal{Y}$, respectively.
It is well known that the eigenvalues of $\linsysMat_0$ are given by   
\begin{equation}\nonumber
\lambda_j(\linsysMat_0)={ 2-2\cos\left(\frac{j 2\pi}{\NI}\right),\quad 
j =  -\left\lfloor \frac{N_{\rm N}-1}{2}\right\rfloor, \dots, -1,0, 1, ..., \left\lfloor \frac{N_{\rm N}}{2}\right\rfloor .}
\end{equation}
For a small argument of the cosine function, we can use the Taylor expansion such that $\lambda_j(\linsysMat_0)=2-2(1-({j2\pi/\NI})^2 /2 )$.
Therefore, the ${2\lfloor {\NI}^\alpha/2\rfloor}+1$ smallest eigenvalues, $j=-\lfloor {\NI}^\alpha/2\rfloor,\dots,-1,0,1,\dots,\lfloor {\NI}^\alpha/2\rfloor$,  are of the same order as $j^2h^2$, where $h = 1/\NI$ denotes the mesh size. 

Let $V_{\alpha}$ contain the  eigenvectors of $\linsysMat_0$ 
(that form a basis of $\mathbb{R}^{\NI}$) except of 
that corresponding to the ${2\lfloor {\NI}^\alpha/2\rfloor}+1$ smallest eigenvalues of $\linsysMat_0$. Then ${\rm dim}(V_{n})={\NI}-{2\lfloor {\NI}^\alpha/2\rfloor}-1$, 
and for $\MB{v}\in V_{\alpha}$ we have
\begin{equation}\nonumber
\frac{\MB{v}^T\linsysMat_{\Jacobi}\MB{v}}
{\MB{v}^T\linsysMat_0\MB{v}}=
\frac{\MB{v}^T(  \linsysMat_0+h^2\MB{T}+
O(h^3)\MB{R}  )\MB{v}}
{\MB{v}^T\linsysMat_0\MB{v}}=1+h^2
\frac{\MB{v}^T  \MB{T}\MB{v}}
{\MB{v}^T\linsysMat_0\MB{v}}
+O(h^3)
\frac{\MB{v}^T  \MB{R}\MB{v}}
{\MB{v}^T\linsysMat_0\MB{v}}=1+O(h^{2\alpha})+O(h),
\end{equation}
where we used the fact that after excluding ${2\lfloor {\NI}^\alpha/2\rfloor}+1$ 
(i.e.~approximatelly $N_{\rm N}^\alpha$) eigenvalues of $\linsysMat_0$,  the smallest remaining eigenvalue of  $\linsysMat_0$ is of order $\frac{{\NI}^{2\alpha}}{\NI^2}=\frac{h^2}{{h}^{2\alpha}}$,
which concludes the proof.

\end{pf}

 \begin{rmk}\label{rmk1}  Let us choose, for example, $\NI=1000$ and $\alpha=1/3$. Then,  according to Lemma~\ref{lemma1}, when we exclude $\NI^{1/3}=10$ largest eigenvalues of the Green-Jacobi preconditioned matrix. The remaining eigenvalues are clustered around unity with a radius of $h^{2/3}\approx 1/100$. This leads to a condition number of approximately$(1+1/100)/(1-1/100)\approx 101/99$. However, minimum of $c$ and magnitudes of $c'$ and $c''$ influence the condition number as well.  The clustering of eigenvalues can accelerate CG convergence, since CG convergence depends on the distribution of eigenvalues rather than just the condition number.
 \end{rmk}

\end{document}